\documentclass[11pt]{article}
\usepackage{amsfonts}
\usepackage{latexsym,amsmath}
\usepackage{amssymb,array}
\usepackage[sort&compress,round,comma]{natbib}
\usepackage{setspace}
\usepackage{enumitem}
\usepackage{xcolor}
\setlist{noitemsep}
\makeatletter \oddsidemargin 0in \evensidemargin 0in \textwidth 16cm
\usepackage{graphicx}
\begin{document}
    \newcommand{\ov}{\overline}
    \newcommand{\om}{\omega}
    \newcommand{\ga}{\gamma}
    \newcommand{\cd}{\circledast}
    \newtheorem{thm}{Theorem}[section]
    \newtheorem{remark}{Remark}[section]
    \newtheorem{counterexample}{Counterexample}[section]
    \newtheorem{coro}{Corollary}[section]
    \newtheorem{propo}{Proposition}[section]
    \newtheorem{definition}{Definition}[section]
    \newtheorem{example}{Example}[section]
    \newtheorem{lem}{Lemma}[section]
    \numberwithin{equation}{section}
    \date{}
\title{Reliability study of series and parallel systems of heterogeneous component lifetimes following proportional odds model}
\author{Pradip Kundu$^1$, Nil Kamal Hazra$^2$ and Asok K. Nanda$^3$\footnote
    {Corresponding author, e-mail:
    asok.k.nanda@gmail.com}\\
    $^1$ Decision Science and Operations Management, Birla Global University,\\ Bhubaneswar, Odisha-751003, India\\
$^2$ Department of Mathematics, IIT Jodhpur, Karwar-342037, India\\
$^3$ Department of Mathematics and Statistics, IISER Kolkata, India}
\maketitle
\begin{abstract}
In this paper, we investigate various stochastic orderings for
series and parallel systems with independent and heterogeneous
components having lifetimes following the proportional odds model.
We also investigate comparisons between system with heterogeneous
components and that with homogeneous components. This paper also
studies relative ageing orders for two systems in the framework of
components having lifetimes following the proportional odds model.
\end{abstract}
\section{Introduction} There is an extensive literature on different stochastic orderings
among order statistics where the observations come from a different
family of distributions. Some of these contributions are due to
\cite{bala}, \cite{bon}, \cite{dyks}, \cite{fangl},
\cite{fang1,fang2}, \cite{gupta}, \cite{kaya}, \cite{khale1},
\cite{khale2}, \cite{koch1, koch2}, \cite{kund}, \cite{li},
\cite{misra}, \cite{nada}, \cite{patr}, \cite{zhao2,zhao3},
\cite{zhao4}, \cite{hkfn,hkfn2}, and the references therein. A
one-to-one correspondence between an order statistic and the
lifetime of a $k$-out-of-$n$ system is well known. A
$k$-out-of-$n:G$ system (generally called $k$-out-of-$n$ system) is
a system consisting of $n$ components which survives as long as at
least $k$ of the $n$ components survive. Let $X_{k:n}$ be the $k$th
smallest order statistic corresponding to the random variables
$X_{1},X_{2},...,X_{n}$, $k=1,2,...,n$. Then the lifetime of an
$(n-k+1)$-out-of-$n:G$ system corresponds to the order statistic
$X_{k:n}$. So, $X_{n-k+1:n}$ represents the lifetime of a
$k$-out-of-$n:G$ system. In particular, $X_{1:n}$ and $X_{n:n}$
represent lifetimes of the series and the parallel systems,
respectively.\par The proportional odds (PO) model introduced by
\cite{benn} is a very important model in survival analysis context,
mainly for its property of convergent hazard functions. The PO
model, as discussed by \cite{benn} and later by \cite{kirmani},
guarantees that the ratio of hazard rates converges to unity as time
tends to infinity. This is in contrast to the proportional hazards
model where the ratio of the hazard rates remains constant with
time. The convergence property of hazard functions makes the PO
model reasonable in many practical applications as discussed in
\cite{benn}, \cite{kirmani} and \cite{ross}. They have also noticed
that assumption of constant hazard ratio is unreasonable in many
practical cases. For more applications of PO model one may refer to
\cite{coll}, \cite{dinse}, \cite{kirmani}, \cite{pett} and the
references therein.\par Let $X$ and $Y$ be two random variables with
distribution functions $F(\cdot)$, $G(\cdot)$, survival functions
$\bar{F}(\cdot)$, $\bar{G}(\cdot)$, probability density functions
$f(\cdot)$, $g(\cdot)$ and hazard rate functions
$r_{X}(\cdot)=f(\cdot)/\bar{F}(\cdot)$,
$r_{Y}(\cdot)=g(\cdot)/\bar{G}(\cdot)$ respectively. Let the odds
functions of $X$ and $Y$ be defined respectively by
$\theta_{X}(t)=\bar{F}(t)/F(t)$ and $\theta_{Y}(t)=\bar{G}(t)/G(t)$.
The random variables $X$ and $Y$ are said to satisfy PO model with
proportionality constant $\alpha$ if $\theta_{Y}(t)=\alpha
\theta_{X}(t)$, for all $t$, where defined. It is observed that, in
terms of survival functions, the PO model can be represented as
\begin{equation}\label{poalt}\bar{G}(t)=\frac{\alpha
\bar{F}(t)}{1-\bar{\alpha}\bar{F}(t)},\end{equation} where
$\bar{\alpha}=1-\alpha$. From the above representation we have
$$\frac{r_{Y}(t)}{r_{X}(t)}=\frac{1}{1-\bar{\alpha}\bar{F}(t)}=\frac{G(t)}{F(t)},$$
so that the hazard ratio is increasing (resp. decreasing) for
$\alpha>1$ (resp. $\alpha<1$) and it converges to unity as $t$ tends
to $\infty$. Also the model (\ref{poalt}), with $0<\alpha <\infty$,
gives a method of introducing new parameter $\alpha$ to a family of
distributions for obtaining more flexible new family of
distributions as discussed in \cite{marsh1}. The family of
distributions so obtained is also known as Marshall-Olkin family of
distributions or Marshall-Olkin extended distributions (for details,
see \cite{marsh1,marsh2} and \cite{cord1} among others).\par
Stochastic comparison of different systems with components following
proportional hazard rates (PHR) model has been discussed by
\cite{dyks}, \cite{khale1}, \cite{koch1,koch2}, and \cite{li} among
others. However, not much work have been done on stochastic
comparison of systems with components following PO model. In this
paper, we investigate stochastic comparisons of series and parallel
systems with heterogeneous components having lifetimes following the
PO model. We also obtain some stochastic comparison results between
a system with heterogeneous components and that with homogeneous
ones. The comparisons are made with respect to the usual stochastic
ordering, the hazard rate ordering, the reversed hazard rate
ordering, the likelihood ratio ordering, and the relative ageing
orderings.\par Throughout the paper, by $a\stackrel{sign}{=} b$ we
mean that $a$ and $b$ have the same sign and by $a\stackrel{def}{=}
b$ we mean that $a$ is defined as $b$. We also write
$\mathbb{R}=(-\infty,\infty)$.
\section{Definitions and preliminaries}\label{prem}
Majorization is a preorder on vectors of real numbers. Let
$I\subseteq \mathbb{R}$ denote a subset of the real line. Further
let, for any vector $\mbox{\boldmath $x$}=(x_1,x_2,...,x_n)\in
\mathbb{R}^n$, $x_{(1)}\leq x_{(2)}\leq ...\leq x_{(n)}$ denote the
increasing arrangement of $x_1,x_2,...,x_n$. Below we give a couple
of definitions to be used throughout the paper.
\begin{definition}
Let $\mbox{\boldmath $x$}=(x_1,x_2,\dots,x_n)\in I^n$ and
$\mbox{\boldmath $y$}=(y_1,y_2,\dots,y_n)\in I^n$. The vector
$\mbox{\boldmath $x$} $ is said to
\begin{enumerate}
\item[(i)] majorize the vector $\mbox{\boldmath $y$}$ (written as $\mbox{\boldmath $x$}\stackrel{m}{\succeq}\mbox{\boldmath $y$}$)
if (cf. \citealp{marsh3})
\begin{equation*}
\sum_{i=1}^j x_{(i)}\le\sum_{i=1}^j y_{(i)},~\text{for
all}\;j=1,\;2,\;\ldots, n-1,\;\;and \;\;\sum_{i=1}^n
x_{(i)}=\sum_{i=1}^n y_{(i)}.
\end{equation*}
\item [(ii)] weakly supermajorize the vector $\mbox{\boldmath $y$}$
 (written as $\mbox{\boldmath $x$}\stackrel{ w}{\succeq} \mbox{\boldmath $y$}$)
 if (cf. \citealp{marsh3})
 \begin{eqnarray*}
  \sum\limits_{i=1}^j x_{(i)}\leq \sum\limits_{i=1}^j y_{(i)},~\text{for all}\;j=1,2,\dots,n.
 \end{eqnarray*}
 \item [(iii)] weakly submajorize the vector $\mbox{\boldmath $y$}$
 (written as $\mbox{\boldmath $x$}\;\succeq_{w} \;\mbox{\boldmath $y$}$)
 if (cf. \citealp{marsh3})
 \begin{eqnarray*}
  \sum\limits_{i=j}^n x_{(i)}\geq \sum\limits_{i=j}^n y_{(i)},~ \text{for all}\;j=1,2,\dots,n.
 \end{eqnarray*}
 \item [(iv)] be $p$-larger than the vector $\mbox{\boldmath $y$}$
 (written as $\mbox{\boldmath $x$}\stackrel{ p}{\succeq} \mbox{\boldmath $y$}$)
 if (cf. \citealp{bon1})
 \begin{eqnarray*}
  \prod\limits_{i=1}^j x_{(i)}\leq \prod\limits_{i=1}^j y_{(i)},~ \text{for all}\;j=1,2,\dots,n.
 \end{eqnarray*}
 \item [(v)] reciprocally majorize the vector $\mbox{\boldmath $y$}$
 (written as $\mbox{\boldmath $x$}\stackrel{ rm}\succeq \mbox{\boldmath $y$}$)
 if (cf. \citealp{zhao1})
 \begin{eqnarray*}
  \sum\limits_{i=1}^j \frac{1}{x_{(i)}}\geq \sum\limits_{i=1}^j \frac{1}{y_{(i)}},~ \text{for all}\;j=1,2,\dots,n.
 \end{eqnarray*}
\end{enumerate}
\end{definition}
It can be seen that
$$\mbox{\boldmath x}\stackrel{m}{\succeq}\mbox{\boldmath
y}\Rightarrow\mbox{\boldmath x}\stackrel{ w}{\succeq}
\mathbf{y}\Rightarrow\mbox{\boldmath x}\stackrel{ p}{\succeq}
\mathbf{y}\Rightarrow\mbox{\boldmath x}\stackrel{ rm}{\succeq}
\mathbf{y}.$$
\begin{remark}\label{remark} Definition 2.1(i) can equivalently be written as
\begin{equation*} \mbox{\boldmath x}\stackrel{m}{\succeq}\mbox{\boldmath
y}~if~ \sum_{i=1}^j x_{[i]}\ge\sum_{i=1}^j y_{[i]},~\text{for
all}\;j=1,\;2,\;\ldots, n-1,\;\;and \;\;\sum_{i=1}^n
x_{[i]}=\sum_{i=1}^n y_{[i]},
\end{equation*} where $x_{[1]}\geq x_{[2]}\geq \cdots \geq x_{[n]}$ is a decreasing arrangement of $x_1,x_2,\cdots,x_n$. \end{remark}
\begin{definition} {\normalfont(\citealp{marsh3})} A function $\phi:I^n \rightarrow \mathbb{R}$ is
said to be Schur-convex (resp. Schur-concave) on $I^n$ if
$$\boldsymbol x\stackrel{m}\succeq \boldsymbol y~\Rightarrow \phi(\boldsymbol x)\geq(\text{resp.} \leq)\phi(\boldsymbol y).$$
\end{definition}
Below we give some definitions of stochastic orders.
\begin{definition}
Let $X$ and $Y$ be two absolutely continuous nonnegative random
variables with cumulative distribution functions $F(\cdot)$,
$G(\cdot)$, survival functions $\bar{F}(\cdot)$, $\bar{G}(\cdot)$,
probability density functions $f(\cdot)$, $g(\cdot)$, hazard rate
functions $r_1(\cdot)$, $r_2(\cdot)$, and the reversed failure
(hazard) rate functions $\tilde r_1(\cdot)$ and $\tilde r_2(\cdot)$,
respectively.
\begin{enumerate}
\item $X$ is said to be smaller than $Y$ in the (cf. \citealp
{shaked})
\begin{enumerate}[label=(\roman*)]
\item usual stochastic order (denoted as $X
\leq_{st} Y$) if $\bar{F}(t)\leq \bar{G}(t)$ for all $t$;
\item failure (hazard) rate order (denoted as $X \leq_{hr} Y$) if
$\bar{G}(t)/\bar{F}(t)$ is increasing in $t\geq 0$, or equivalently
if $r_1(t)\geq r_2(t)$ for all $t\geq 0$;
\item reversed failure (hazard) rate order (denoted as $X \leq_{rhr} Y$) if
$G(t)/F(t)$ is increasing in $t>0$, or equivalently if $\tilde
r_1(t)\leq \tilde r_2(t)$ for all $t>0$;
\item likelihood ratio order (denoted as $X \leq_{lr} Y$) if $f(x)/g(x)$ decreases in $x$ over the union of the supports of $X$ and
$Y$.
\end{enumerate}
\item  $X$ is said to age faster than $Y$ in terms of the
\begin{enumerate}[label=(\roman*)]
\item hazard rate (denoted as $X\lesssim_{hr}Y$), if
$r_1(t)/r_2(t)$ is increasing in $t>0$ (cf. \citealp{seng});
\item reversed hazard rate (denoted as $X\lesssim_{rhr}Y$), if
$\tilde r_2(t)/\tilde r_1(t)$ is increasing in $t>0$ (cf.
\citealp{reza}).$\hfill\Box$
\end{enumerate}
\end{enumerate}
\end{definition}
The following notations are used throughout the paper.
\begin{enumerate}[label=(\roman*)]
\item $\mathcal{D}=\{(x_1,x_2,...,x_n)\in \mathbb{R}^n :x_1\geq x_2\geq \cdots \geq x_n
\}$.
\item $\mathcal{D_{+}}=\{(x_1,x_2,...,x_n)\in \mathbb{R}^n :x_1\geq x_2\geq \cdots \geq
x_n>0 \}$.
\item $\mathcal{E}=\{(x_1,x_2,...,x_n)\in \mathbb{R}^n :x_1\leq x_2\leq \cdots \leq x_n
\}$.
\item $\mathcal{E_{+}}=\{(x_1,x_2,...,x_n)\in \mathbb{R}^n :0<x_1\leq x_2\leq \cdots \leq
x_n\}$.
\end{enumerate}
Before we start, we mention below, for completeness, a few lemmas to
be used in the sequel. Below we take ${\mbox{\boldmath
$z$}}=(z_1,z_2,...,z_n)$ and $\varphi_{(k)}({\mbox{\boldmath
$z$}})=\partial\varphi({\mbox{\boldmath $z$}})/\partial z_k$, the
partial derivative of $\varphi$ with respect to its $k$th argument.
\begin{lem}\label{le1} {\normalfont(\citealp{marsh3})}
 Let $\varphi:\mathcal{D}\rightarrow \mathbb{R}$ be a function, continuously differentiable on the interior of $\mathcal{D}$. Then, for ${\mbox{\boldmath $x$}},{\mbox{\boldmath $y$}}\in \mathcal{D}$,
 \begin{eqnarray*}
  {\mbox{\boldmath $x$}}\stackrel{m}\succeq{\mbox{\boldmath $y$}}\;\implies\;\varphi({\mbox{\boldmath $x$}})\geq (\text{resp.}\;\leq)\;\varphi({\mbox{\boldmath $y$}})
 \end{eqnarray*}
if, and only if,
$$\varphi_{(k)}({\mbox{\boldmath $z$}})\;\text{is decreasing (resp. increasing) in}\;k=1,2,\dots,n.$$
\end{lem}
%%%%%%%%%%%%%%%%%%%%%% Lemma 3.2 %%%%%%%%%%%%%%%%%%%%%%%%%%%%%%%%%
\begin{lem}\label{le2} {\normalfont(\citealp{marsh3})}
 Let $\varphi:\mathcal{E}\rightarrow \mathbb{R}$ be a function, continuously differentiable on the interior of $\mathcal{E}$.
 Then, for ${\mbox{\boldmath $x$}},{\mbox{\boldmath $y$}}\in \mathcal{E}$,
 \begin{eqnarray*}
  {\mbox{\boldmath $x$}}\stackrel{m}\succeq{\mbox{\boldmath $y$}}\;\implies\;\varphi({\mbox{\boldmath $x$}})\geq (\text{resp.}\;\leq)\;\varphi({\mbox{\boldmath $y$}})
 \end{eqnarray*}
if, and only if,
$$\varphi_{(k)}({\mbox{\boldmath $z$}})\;\text{is increasing (resp. decreasing) in}\;k=1,2,\dots,n.$$
\end{lem}
%%%%%%%%%%%%%%%%%%%%%%%%  Lemma 3.3  %%%%%%%%%%%%%%%%%%%%%%%%%%%%%%
\begin{lem}\label{le2a} {\normalfont(\citealp{marsh3})}
 Let $I \subseteq \mathbb{R}$ be an open interval and let $\varphi: I^n\rightarrow \mathbb{R}$ be continuously differentiable. Necessary
and sufficient conditions for $\varphi$ to be Schur-convex (resp.
Schur-concave) on $I^n$ are $\varphi$ is symmetric on $I^n$, and for
all $i\neq j$,
$$(z_i-z_j)\left(\varphi_{(i)}({\mbox{\boldmath
$z$}})-\varphi_{(j)}({\mbox{\boldmath $z$}})\right)\geq
(\text{resp.}\leq)\;0~\text{for all}~{\mbox{\boldmath $z$}}\in I^n.$$
 \end {lem}
%%%%%%%%%%%%%%%%%%%%%%%%  Lemma 3.3  %%%%%%%%%%%%%%%%%%%%%%%%%%%%%%
 \begin{lem}\label{le3} {\normalfont(\citealp{marsh3})}
 Let $S \subseteq \mathbb{R}^n$. Further, let $\varphi: S\rightarrow \mathbb{R}$ be a function. Then, for ${\mbox{\boldmath $x$}},{\mbox{\boldmath $y$}}\in S,$
 $${\mbox{\boldmath $x$}}\succeq_{w} {\mbox{\boldmath $y$}}\; \implies\;\varphi({\mbox{\boldmath $x$}})\geq (\text{resp. }\leq)\; \varphi({\mbox{\boldmath $y$}})$$
 if, and only if, $\varphi$ is both increasing (resp. decreasing) and Schur-convex (resp. Schur-concave) on $S$.
 Similarly,
$${\mbox{\boldmath $x$}}\stackrel{w}\succeq {\mbox{\boldmath $y$}}\; \implies\;\varphi({\mbox{\boldmath $x$}})\geq (\text{resp. }\leq)\; \varphi({\mbox{\boldmath $y$}})$$
if, and only if, $\varphi$ is both decreasing (resp. increasing) and
Schur-convex (resp. Schur-concave) on $S$. $\hfill\Box$
 \end {lem}
%%%%%%%%%%%%%%%%%%%%%%%  Lemma 3.4  %%%%%%%%%%%%%%%%%%%%%%%%%%%%%%%%
\begin{lem}\label{le4} {\normalfont(\citealp{khale}; \citealp{kund})}
  Let $\varphi :(0,\infty)^n\rightarrow \mathbb{R}$ be a function.
 Then,
 $${\mbox{\boldmath $x$}}\stackrel{\rm p}\succeq {\mbox{\boldmath $y$}}\implies \varphi({\mbox{\boldmath $x$}})\geq(\text{resp.}\leq) \;\varphi({\mbox{\boldmath $y$}})$$
if, and only if, the following two conditions hold:
\begin{enumerate}
 \item[(i)]$\varphi(e^{a_1},\dots,e^{a_n})$ is Schur-convex (resp. Schur-concave) in $(a_1,\dots,a_n)$,
\item[(ii)]$\varphi(e^{a_1},\dots,e^{a_n})$ is decreasing (resp. increasing) in each $a_i,$ for $i=1,\dots,n,$
\end{enumerate}
where $a_i=\ln x_i$, for $i=1,\dots,n.$ $\hfill\Box$
\end{lem}
%%%%%%%%%%%%%%%%%%%%   Lemma- 3.5  %%%%%%%%%%%%%%%%%%%%%%%%%%%%%%%%
Following lemma is adapted from \cite{bon} (see also
\citealp{gupta}).
\begin{lem}\label{lemgm} Let $\phi:(0,\infty)^{n}\rightarrow (0,\infty)$ be a symmetrical and continuously
differentiable mapping. If, for $\textbf{x}=(x_1,x_2,...,x_n)\in
(0,\infty)^{n}$ with $x_p=\min_{1\leq i\leq n}x_i$ and
$x_q=\max_{1\leq i\leq n}x_i$, we have
$$(x_p-x_q)\left(\frac{1}{\prod_{i\neq p} x_i}\frac{\partial\phi}{\partial
x_p}-\frac{1}{\prod_{i\neq q} x_i}\frac{\partial\phi}{\partial
x_q}\right)<(>) 0,$$ for $x_p\neq x_q$, then
$$\phi(x_1,x_2,...,x_n)\leq (\geq) \phi(x,x,...,x),$$ where
$x=\sqrt[n]{x_1 x_2\cdots x_n}$.\end{lem}
\section{Series systems with component lifetimes following PO model}
In this section we compare the lifetimes of two series systems, each
of the heterogeneous components having lifetimes following the PO
model, with respect to some stochastic orders. We also compare
lifetimes of two series systems, one comprising of heterogeneous
components and another comprising of homogeneous components.\par
Throughout the paper we consider $X=(X_1,X_2,...,X_n)$ and
$Y=(Y_1,Y_2,...,Y_n)$ as two sets of independent random variables.
Let both $X$ and $Y$ follow the PO model, denoted as $X\sim
PO(\bar{F},\boldsymbol\lambda)$ and $Y\sim
PO(\bar{F},\boldsymbol\mu)$, where $\bar{F}$ is the baseline
survival function,
$\boldsymbol\lambda=(\lambda_1,\lambda_2,...,\lambda_n)$ and
$\boldsymbol\mu=(\mu_1,\mu_2,...,\mu_n)$ with $\lambda_{i}>0$ and
$\mu_i>0$, for all $i=1,2,...,n$. We have the survival functions of
$X_{1:n}$ and $Y_{1:n}$, respectively, as
$$\bar{F}_{X_{1:n}}(x)=\prod_{i=1}^{n} \bar{F}_{X_i}(x)=\prod_{i=1}^{n}\frac{\lambda_i
\bar{F}(x)}{1-\bar{\lambda}_i\bar{F}(x)},$$ and
$$\bar{F}_{Y_{1:n}}(x)=\prod_{i=1}^{n} \bar{F}_{Y_i}(x)=\prod_{i=1}^{n}\frac{\mu_i
\bar{F}(x)}{1-\bar{\mu}_i\bar{F}(x)},$$ where
$\bar{\lambda}_i=1-\lambda_i$ and $\bar{\mu}_i=1-\mu_i$, for
$i=1,2,...,n$.

The hazard rate functions of $X_{1:n}$ and
$Y_{1:n}$ are, respectively, obtained as
$$r_{X_{1:n}}(x)=\sum_{i=1}^{n} r_{X_i}(x)=\sum_{i=1}^n
\frac{r(x)}{1-\bar{\lambda}_i\bar{F}(x)},$$ and
$$r_{Y_{1:n}}(x)=\sum_{i=1}^{n} r_{Y_i}(x)=\sum_{i=1}^n
\frac{r(x)}{1-\bar{\mu}_i\bar{F}(x)}.$$ If $X\sim PO(\bar{F},\lambda
\boldsymbol 1)$, where $\boldsymbol 1=(1,1,...,1)$, $\lambda>0$,
then the survival function and the hazard rate function of $X_{1:n}$
are given respectively by
$$\bar{F}_{X_{1:n}}(x)=\frac{\lambda^n
\bar{F}^n(x)}{(1-\bar{\lambda}\bar{F}(x))^n},$$ and
$$r_{X_{1:n}}(x)=\frac{n r(x)}{1-\bar{\lambda}\bar{F}(x)},$$
where $\bar{\lambda}=1-\lambda$.

Suppose each of the two series systems is formed out of $n$
heterogeneous components where the component lifetimes follow the PO
model. The following theorem compares the lifetimes of two such
series systems.
\begin{thm}\label{thst}Suppose the lifetime vectors $X\sim PO(\bar{F},\boldsymbol\lambda)$ and $Y\sim
PO(\bar{F},\boldsymbol\mu)$. Then
$$\boldsymbol\lambda \stackrel{p}\succeq\boldsymbol\mu~\text{implies}~ X_{1:n}\leq_{st}Y_{1:n}.$$
\end{thm}
\textbf{Proof:} Write $a_i=\ln \lambda_i$, $i=1,2,...,n$. Then
\begin{eqnarray*}\bar{F}_{X_{1:n}}(x)&=&\prod_{i=1}^{n}\frac{e^{a_i}
\bar{F}(x)}{1-(1-e^{a_i})\bar{F}(x)}\\&=&\phi(e^{a_1},e^{a_2},...,e^{a_n}),~\text{(say)}.\end{eqnarray*}
Note that $\phi(e^{a_1},e^{a_2},...,e^{a_n})$ is symmetric with
respect to $(a_1,a_2,...,a_n)\in \mathbb{R}^n$. Now,
$$\frac{\partial\phi}{\partial
a_i}=\frac{1-\bar{F}(x)}{1-(1-e^{a_i})\bar{F}(x)}\phi(e^{a_1},e^{a_2},...,e^{a_n}),$$
so that $\phi(e^{a_1},e^{a_2},...,e^{a_n})$ is increasing in each
$a_i$, for $i=1,2,...,n$. Now, for $1\leq i\leq j\leq n$,
\begin{eqnarray*}(a_i- a_j)\left(\frac{\partial\phi}{\partial
a_i}-\frac{\partial\phi}{\partial
a_j}\right)&=&\frac{(a_i- a_j)(e^{a_j}-e^{a_i})\bar{F}(x)(1-\bar{F}(x))}{(1-(1-e^{a_i})\bar{F}(x))(1-(1-e^{a_j})\bar{F}(x))}\phi(e^{a_1},e^{a_2},...,e^{a_n})\\
&\leq& 0.\end{eqnarray*} So, from Lemma \ref{le2a},
$\phi(e^{a_1},e^{a_2},...,e^{a_n})$ is Schur-concave in
$(a_1,a_2,...,a_n)\in \mathbb{R}^n$. Thus, from Lemma \ref{le4}, we
have $\phi(\lambda_1,\lambda_2,...,\lambda_n)\leq \phi
(\mu_1,\mu_2,...,\mu_n)$ whenever $\boldsymbol\lambda
\stackrel{p}\succeq \boldsymbol\mu$. This proves the result.$\hfill\Box$

The following corollary immediately follows from the above theorem
by noting the fact that
$(\lambda_1,\lambda_2,\ldots,\lambda_n)\stackrel{p}\succeq(\underbrace{\lambda,\lambda,\ldots,\lambda}_{n\;terms})$,
where $\lambda\geqslant \left(\prod_{i=1}^n\lambda_i\right)^{1/n}.$
\begin{coro}\label{corst}Suppose that the lifetime vectors $X\sim PO(\bar{F},\boldsymbol\lambda)$ and $Y\sim
PO(\bar{F},\lambda \boldsymbol 1)$. Then, $X_{1:n}\leq_{st}Y_{1:n}$
if $\lambda \geq \sqrt[n]{\lambda_1 \lambda_2 \cdots
\lambda_n}$.$\hfill\Box$
\end{coro}
Since $p$-larger order is stronger than reciprocal majorization
order, one may wonder whether, in Theorem \ref{thst}, $p$-larger
order can be replaced by reciprocal majorization order. The
Counterexample \ref{nsturm} shows that this cannot be done.

Since hazard rate order is stronger than stochastic order, in order
to get a comparison of series systems in terms of hazard rate order,
we need to have some larger dominance than the $p$-larger order
between the parameters of the models. The following theorem gives a
condition under which two series systems formed out of component
lifetimes following the PO models will be ordered in hazard rate
order.
\begin{thm}\label{thhr}Suppose that the lifetime vectors $X\sim PO(\bar{F},\boldsymbol\lambda)$ and $Y\sim
PO(\bar{F},\boldsymbol\mu)$. Then
$$\boldsymbol\lambda \stackrel{w}\succeq \boldsymbol\mu~\text{implies}~ X_{1:n}\leq_{hr}Y_{1:n}.$$ \end{thm}
\textbf{Proof:} We have $$r_{X_{1:n}}(x)=\sum_{i=1}^n
\frac{r(x)}{1-\bar{\lambda}_i\bar{F}(x)},$$ which is symmetric with
respect to $(\lambda_1,\lambda_2,...,\lambda_n)\in \mathbb{R}^n$. Differentiating the above expression with respect to $\lambda_i$ we get
$$\frac{\partial r_{X_{1:n}}(x)}{\partial
\lambda_i}=-\frac{r(x)\bar{F}(x)}{(1-\bar{\lambda}_i\bar{F}(x))^
2}\;<0,$$ which tells that $r_{X_{1:n}}(x)$ is decreasing in each
$\lambda_i$, $i=1,2,...,n.$ For $1\leq i\leq j\leq n$,
\begin{eqnarray*}(\lambda_i-\lambda_j)\left(\frac{\partial r_{X_{1:n}}(x)}{\partial
\lambda_i}-\frac{\partial r_{X_{1:n}}(x)}{\partial
\lambda_j}\right)&=&(\lambda_i-\lambda_j)r(x)\bar{F}(x)\left[\frac{1}{(1-\bar{\lambda_j}\bar{F}(x))^
2}-\frac{1}{(1-\bar{\lambda}_i\bar{F}(x))^2}\right]\\&\stackrel{sign}{=}&(\lambda_i-\lambda_j)\left[(1-\bar{\lambda}_i\bar{F}(x))^2-(1-\bar{\lambda_j}\bar{F}(x))^
2\right]\\&\geq& 0.
\end{eqnarray*}
So, from Lemma \ref{le2a}, it follows that $r_{X_{1:n}}(x)$ is
Schur-convex in
$\boldsymbol\lambda=(\lambda_1,\lambda_2,...,\lambda_n)\in
\mathbb{R}^n$. Thus, by Lemma \ref{le3}, we have $r_{X_{1:n}}(x)\geq
r_{Y_{1:n}}(x)$ whenever $\boldsymbol\lambda
\stackrel{w}\succeq\boldsymbol\mu$.$\hfill\Box$\\
Since
$(\lambda_1,\lambda_2,\ldots,\lambda_n)\stackrel{w}\succeq(\underbrace{\lambda,\lambda,\ldots,\lambda}_{n\;terms})$,
for $\lambda\geqslant\frac{1}{n}\sum_{i=1}^n\lambda_i$, the
following corollary immediately follows from the above theorem.
\begin{coro}\label{corhr}Suppose lifetime vectors $X\sim PO(\bar{F},\boldsymbol\lambda)$ and $Y\sim
PO(\bar{F},\lambda \boldsymbol 1)$. Then, $X_{1:n}\leq_{hr}Y_{1:n}$
if $\lambda \geq \frac{1}{n}\sum_{i=1}^n \lambda_i$.
$\hfill\Box$
\end{coro}

Since weakly supermajorization order is stronger than $p$-larger
order, one may wonder whether weakly supermajorization order in
Theorem \ref{thhr} can be replaced by $p$-larger order. The
Counterexample \ref{nhrup} shows that this cannot be
done.$\hfill\Box$

Let $X_1,X_2,\ldots,X_p$ have a common distribution $F$ and let
$X_{p+1},X_{p+2},\ldots,X_n$ have a common distribution $G$, for
$p=1,2,\ldots,n-1$. The distribution $F$ is called the original
distribution whereas the distribution $G$ is called the outlier
distribution. This type of model is known as outlier model. For
$p=n-1$, the model is known as a single-outlier model whereas, for
$p=1,2,\ldots,n-2$, the model is called multiple-outlier model.
Below we study the relative ageing of two series systems with
heterogeneous components in terms of the hazard rate in the case of
multiple-outlier model.
\begin{thm}\label{thag}Let both $X$ and $Y$ follow the multiple-outlier PO
model with $X_i\sim PO(\bar{F},\lambda_1)$, $Y_i\sim
PO(\bar{F},\mu_1)$, for $i=1,2,...,n_1$, $X_j\sim
PO(\bar{F},\lambda_2)$, $Y_j\sim PO(\bar{F},\mu_2)$, for
$j=n_1+1,n_1+2,...,n_1+n_2(=n)$. Then
$$(\underbrace{\lambda_1,\lambda_1,...,\lambda_1}_{n_1\;terms},\underbrace{\lambda_2,\lambda_2,...,\lambda_2}_{n_2\;terms})\stackrel{m}\succeq
(\underbrace{\mu_1,\mu_1,...,\mu_1}_{n_1\;terms},\underbrace{\mu_2,\mu_2,...,\mu_2}_{n_2\;terms})\Rightarrow
X_{1:n}\gtrsim_{hr}Y_{1:n},$$ provided $\{(\lambda_1,\lambda_2)\in
\mathcal{E_{+}},(\mu_1,\mu_2)\in \mathcal{E_{+}}\}$ or
$\{(\lambda_1,\lambda_2)\in \mathcal{D_{+}},(\mu_1,\mu_2)\in
\mathcal{D_{+}}\}$.
\end{thm}
\textbf{Proof:} First it should be noted that
$$(\underbrace{\bar{\lambda}_1,\bar{\lambda}_1,...,\bar{\lambda}_1}_{n_1\;terms},\underbrace{\bar{\lambda}_2,\bar{\lambda}_2,...,\bar{\lambda}_2}_{n_2\;terms})\stackrel{m}\succeq(\underbrace{\bar{\mu}_1,\bar{\mu}_1,...,\bar{\mu}_1}_{n_1\;terms},\underbrace{\bar{\mu}_2,\bar{\mu}_2,...,\bar{\mu}_2}_{n_2\;terms})$$
is equivalent to
$$(\underbrace{\lambda_1,\lambda_1,...,\lambda_1}_{n_1\;terms},\underbrace{\lambda_2,\lambda_2,...,\lambda_2}_{n_2\;terms})\stackrel{m}\succeq
(\underbrace{\mu_1,\mu_1,...,\mu_1}_{n_1\;terms},\underbrace{\mu_2,\mu_2,...,\mu_2}_{n_2\;terms}),$$
which follows from Remark 2.1 and Definition 2.1(i). Here we write $\bar{\lambda}_i=1-\lambda_i$ and $\bar{\mu}_i=1-\mu_i$, for $i=1,2.$\\
We denote
$$\mathcal{A}=\{\boldsymbol\xi=(\xi_1,\xi_2,...,\xi_n):\xi_i=\lambda_1~\text{for}~ 1\leq i\leq n_1~\text{and}~\xi_j=\lambda_2,~\text{for}~n_1+1\leq j\leq n\}$$
and
$$\mathcal{B}=\{\boldsymbol\eta=(\eta_1,\eta_2,...,\eta_n):\eta_i=\mu_1~\text{for}~ 1\leq i\leq n_1~\text{and}~\eta_j=\mu_2,~\text{for}~n_1+1\leq j\leq n\}.$$
We have to show that, under the given majorization order,
\begin{equation}\label{agfh}\frac{r_{X_{1:n}}(x)}{r_{Y_{1:n}}(x)}=\frac{\sum_{i=1}^n
\frac{1}{1-\bar{\xi}_i\bar{F}(x)}}{\sum_{i=1}^n
\frac{1}{1-\bar{\eta}_i\bar{F}(x)}}\;{\rm is\; decreasing\; in}\;x>0,
\end{equation}
for all $\boldsymbol\xi\in \mathcal{A}$, $\boldsymbol\eta\in
\mathcal{B}$, where $\bar{\xi}_i=1-\xi_i$ and
$\bar{\eta}_i=1-\eta_i$, for $i=1,2,\ldots,n$, which is equivalent
to show that
$$\frac{\sum_{i=1}^n
\frac{\bar{\xi}_i}{(1-\bar{\xi}_i\bar{F}(x))^2}}{\sum_{i=1}^n
\frac{1}{1-\bar{\xi}_i\bar{F}(x))}}\geq \frac{\sum_{i=1}^n
\frac{\bar{\eta}_i}{(1-\bar{\eta}_i\bar{F}(x))^2}}{\sum_{i=1}^n
\frac{1}{1-\bar{\eta}_i\bar{F}(x)}}.$$ Furthermore, to prove this,
it suffices to show that (according to Definition 2.2), for all
$\boldsymbol\xi\in \mathcal{A}$ and $\boldsymbol\eta\in
\mathcal{B}$,
$$\phi(\bar{\xi}_1,\bar{\xi}_2,...,\bar{\xi}_n)\stackrel{ def}=\frac{\sum_{i=1}^n
\frac{\bar{\xi}_i}{(1-\bar{\xi}_i\bar{F}(x))^2}}{\sum_{i=1}^n
\frac{1}{1-\bar{\xi}_i\bar{F}(x)}}$$ is Schur-convex in
$(\bar{\xi}_1,\bar{\xi}_2,...,\bar{\xi}_n)\in \mathcal{A}$. Now, writing
$u(x)=1/(1-x)$ and $v(x)=x/(1-x)$, we have
\begin{eqnarray*}\phi(\bar{\xi}_1,\bar{\xi}_2,...,\bar{\xi}_n)&=&\frac{1}{\bar{F}(x)}\frac{\sum_{i=1}^n u(\bar{\xi}_i\bar{F}(x))v(\bar{\xi}_i\bar{F}(x))}{\sum_{i=1}^n u(\bar{\xi}_i\bar{F}(x))}\\
&=&\frac{1}{\bar{F}(x)}\frac{\sum_{i=1}^{n_1}
u(\bar{\xi}_i\bar{F}(x))v(\bar{\xi}_i\bar{F}(x))+\sum_{i=n_1+1}^{n}
u(\bar{\xi}_i\bar{F}(x))v(\bar{\xi}_i\bar{F}(x))}{\sum_{i=1}^{n_1}
u(\bar{\xi}_i\bar{F}(x))+\sum_{i=n_1+1}^{n}
u(\bar{\xi}_i\bar{F}(x))}.\end{eqnarray*}
 Writing $u'=du/dx$ and $v'=dv/dx$, we have, for $1\leq i\leq n_1$,
{\scriptsize$$\frac{\partial \phi}{\partial \bar{\xi}_i}=\frac{
u(\bar{\lambda}_1\bar{F}(x))v'(\bar{\lambda}_1\bar{F}(x))[n_1
u(\bar{\lambda}_1\bar{F}(x))+n_2 u(\bar{\lambda}_2\bar{F}(x))]+ n_2
u(\bar{\lambda}_2\bar{F}(x))
u'(\bar{\lambda}_1\bar{F}(x))[v(\bar{\lambda}_1\bar{F}(x))-v(\bar{\lambda}_2\bar{F}(x))]}{\left(n_1
u(\bar{\lambda}_1\bar{F}(x))+n_2
u(\bar{\lambda}_2\bar{F}(x))\right)^2},$$} \text{ and, for } $n_1
+1\leq j \leq n$, {\scriptsize$$\frac{\partial \phi}{\partial
\bar{\xi}_j}=\frac{
u(\bar{\lambda}_2\bar{F}(x))v'(\bar{\lambda}_2\bar{F}(x))[n_1
u(\bar{\lambda}_1\bar{F}(x))+n_2 u(\bar{\lambda}_2\bar{F}(x))]+n_1
 u(\bar{\lambda}_1\bar{F}(x))
u'(\bar{\lambda}_2\bar{F}(x))[v(\bar{\lambda}_2\bar{F}(x))-v(\bar{\lambda}_1\bar{F}(x))]}{\left(n_1
u(\bar{\lambda}_1\bar{F}(x))+n_2
u(\bar{\lambda}_2\bar{F}(x))\right)^2}.$$} Now, for $1\leq i,j\leq
n_1$ or $n_1 +1\leq i,j \leq n$, we have
$\displaystyle\frac{\partial \phi}{\partial
\bar{\xi}_i}-\frac{\partial \phi}{\partial \bar{\xi}_j}=0$. Again,
for $1\leq i\leq n_1$ and $n_1 +1\leq j \leq n$, we have
\allowdisplaybreaks{
\begin{eqnarray*}
\frac{\partial \phi}{\partial \bar{\xi}_i}-\frac{\partial
\phi}{\partial \bar{\xi}_j}&\stackrel{sign}{=}& [n_1
u(\bar{\lambda_1}\bar{F}(x))+n_2 u(\bar{\lambda}_2\bar{F}(x))][
u(\bar{\lambda}_1\bar{F}(x))v'(\bar{\lambda}_1\bar{F}(x))-
u(\bar{\lambda}_2\bar{F}(x))v'(\bar{\lambda}_2\bar{F}(x))]\\&&+
[v(\bar{\lambda}_1\bar{F}(x))-v(\bar{\lambda}_2\bar{F}(x))][n_2u(\bar{\lambda}_2\bar{F}(x))
u'(\bar{\lambda}_1\bar{F}(x))+n_1u(\bar{\lambda}_1\bar{F}(x))
u'(\bar{\lambda}_2\bar{F}(x))].
\end{eqnarray*}} Since $v(x)$ and
$u(x)v'(x)$ are both increasing in $x$, $u(x)$ is nonnegative for
all $x\leq 1$ and $u'(x)$ is nonnegative for all $x$, we have, for
$\bar{\lambda}_1\geq(\text{resp.}\leq)~ \bar{\lambda}_2$,
$$\frac{\partial \phi}{\partial \bar{\xi}_i}-\frac{\partial
\phi}{\partial \bar{\xi}_j}\geq(\text{resp.} \leq)~0.$$ So, from
Lemma \ref{le1} and Lemma \ref{le2}, it follows that $\phi$ is
Schur-convex in $(\bar{\xi}_1,\bar{\xi}_2,...,\bar{\xi}_n)\in
\mathcal{A}$. Thus,%, under the set of
%conditions $\{(\bar{\lambda}_1,\bar{\lambda}_2)\in
%\mathcal{D},(\bar{\mu}_1,\bar{\mu}_2)\in \mathcal{D},n_1\geq n_2\}$
%or $\{(\bar{\lambda}_1,\bar{\lambda}_2)\in
%\mathcal{E},(\bar{\mu}_1,\bar{\mu}_2)\in \mathcal{E},n_1\leq n_2\}$,
$$(\underbrace{\bar{\lambda}_1,\bar{\lambda}_1,...,\bar{\lambda}_1}_{n_1\;terms},\underbrace{\bar{\lambda}_2,\bar{\lambda}_2,...,\bar{\lambda}_2}_{n_2\;terms})\stackrel{m}\succeq(\underbrace{\bar{\mu}_1,\bar{\mu}_1,...,\bar{\mu}_1}_{n_1\;terms},\underbrace{\bar{\mu}_2,\bar{\mu}_2,...,\bar{\mu}_2}_{n_2\;terms})
\Rightarrow X_{1:n}\gtrsim_{hr}Y_{1:n},$$ and hence the result is proved.\hfill$\Box$
%%%%%%%%%%%%%%%%%%%%%%%%%%%%%%%%%%%%%%%%%%%%%

By taking $n_1=n_2=1$ in the above theorem, we immediately get the
following corollary.
\begin{coro}Let, for $i=1,2$, the two independent random variables $X_i$ and $Y_i$ follow the PO model with parameters $\lambda_i$ and $\mu_{i}$, respectively. Then
$$(\lambda_1,\lambda_2)\stackrel{m}\succeq(\mu_1,\mu_2)\Rightarrow X_{1:2}\gtrsim_{hr}Y_{1:2}.$$\end{coro}
One may wonder whether the set of sufficient conditions given in
Theorem \ref{thag} is the only possible set of conditions or the
result is possible to be true under a different set of sufficient
conditions. The following theorem answers this in affirmative.
\begin{thm}\label{thag2}Let both $X$ and $Y$ follow the multiple-outlier PO
model with $X_i\sim PO(\bar{F},\lambda_1)$, $Y_i\sim
PO(\bar{F},\mu_1)$, for $i=1,2,...,n_1$, $X_j\sim
PO(\bar{F},\lambda_2)$, $Y_j\sim PO(\bar{F},\mu_2)$, for
$j=n_1+1,n_1+2,...,n_1+n_2(=n)$. Then
$$\max\{\lambda_1,\lambda_2\}\leq \min\{\mu_1,\mu_2\}\Rightarrow
X_{1:n}\gtrsim_{hr}Y_{1:n}.$$
\end{thm}
\textbf{Proof:} Proving $X_{1:n}\gtrsim_{hr}Y_{1:n}$ is equivalent
to showing that
\begin{equation}\label{agfh2}\frac{r_{X_{1:n}}(x)}{r_{Y_{1:n}}(x)}=\frac{
\frac{n_1}{1-\bar{\lambda_1}\bar{F}(x)}+\frac{n_2}{1-\bar{\lambda_2}\bar{F}(x)}}{
\frac{n_1}{1-\bar{\mu_1}\bar{F}(x)}+\frac{n_2}{1-\bar{\mu_2}\bar{F}(x)}}
~\text{is decreasing in}~ x>0.\end{equation} Writing $u(x)=1/(1-x)$
and $v(x)=x/(1-x)$, (\ref{agfh2}) becomes equivalent to the fact
that
\begin{eqnarray*}&&n_1^2
u(\bar{\mu_1}\bar{F}(x))u(\bar{\lambda_1}\bar{F}(x))[v(\bar{\mu_1}\bar{F}(x))-v(\bar{\lambda_1}\bar{F}(x))]+n_1
n_2 u(\bar{\mu_1}\bar{F}(x))u(\bar{\lambda_2}\bar{F}(x))
\\&&[v(\bar{\mu_1}\bar{F}(x))-v(\bar{\lambda_2}\bar{F}(x))]+n_1 n_2
u(\bar{\mu_2}\bar{F}(x))u(\bar{\lambda_1}\bar{F}(x))[v(\bar{\mu_2}\bar{F}(x))-v(\bar{\lambda_1}\bar{F}(x))]\\&&+n_2^2
u(\bar{\mu_2}\bar{F}(x))u(\bar{\lambda_2}\bar{F}(x))[v(\bar{\mu_2}\bar{F}(x))-v(\bar{\lambda_2}\bar{F}(x))]\leq
0.\end{eqnarray*} As both $u(x)$ and $v(x)$ are increasing in $x$, the above inequality holds if the condition
$\max\{\lambda_1,\lambda_2\}\leq \min\{\mu_1,\mu_2\}$
holds. This proves the theorem.
\begin{remark}
From Theorem \ref{thag}, we get that  $X_{1:n}\gtrsim_{hr}Y_{1:n}$
whenever
    \begin{equation}\label{aa}
    (\underbrace{\lambda_1,\lambda_1,...,\lambda_1}_{n_1\;terms},\underbrace{\lambda_2,\lambda_2,...,\lambda_2}_{n_2\;terms})\stackrel{m}\succeq
    (\underbrace{\mu_1,\mu_1,...,\mu_1}_{n_1\;terms},\underbrace{\mu_2,\mu_2,...,\mu_2}_{n_2\;terms})
    \end{equation}
    whereas, from Theorem \ref{thag2}, we have that $X_{1:n}\gtrsim_{hr}Y_{1:n}$ if
    \begin{equation}\label{bb}
    \max\{\lambda_1,\lambda_2\}\leq \min\{\mu_1,\mu_2\}.
    \end{equation}
    From these two theorems one natural question that arises is -- whether
    (\ref{aa}) $\Rightarrow$ (\ref{bb}) or (\ref{bb}) $\Rightarrow$ (\ref{aa}). This is because if (\ref{aa}) $\Rightarrow$ (\ref{bb}) then Theorem \ref{thag} is redundant whereas if (\ref{bb}) $\Rightarrow$ (\ref{aa}), Theorem \ref{thag2} will be redundant. By taking
     $\lambda_1=1$, $\lambda_2=2$, $\mu_1=1.3$, $\mu_2=1.8$, $n_1=2$ and $n_2=3$, it is clear that (\ref{aa}) is satisfied but not (\ref{bb}). Further, (\ref{bb}) cannot imply (\ref{aa}) because if (\ref{bb}) holds then $n_1\lambda_1+n_2\lambda_2=n_1\mu_1+n_2\mu_2$ is never satisfied and hence (\ref{aa}) cannot hold.\hfill$\Box$
\end{remark}
    Looking into Theorem \ref{thag}, one may wonder whether the condition on majorization order can be relaxed to the weakly supermajorization order. Below we answer this question in affirmative. However, for this relaxation, we need to sacrifice the broadness of the model in terms of the parameters. %A result on relative ageing is given next in terms of weakly majorization order.
\begin{thm}\label{thag3}Let both $X$ and $Y$ follow the multiple-outlier PO
model with $X_i\sim PO(\bar{F},\lambda_1)$, $Y_i\sim
PO(\bar{F},\mu_1)$, for $i=1,2,...,n_1$, $X_j\sim PO(\bar{F},\eta)$,
$Y_j\sim PO(\bar{F},\eta)$, for $j=n_1+1,n_1+2,...,n_1+n_2(=n)$.
Then for $\lambda_1\leq \min\{\eta,\mu_1\}$,
$$(\underbrace{\lambda_1,\lambda_1,...,\lambda_1}_{n_1\;terms},\underbrace{\eta ,\eta ,...,\eta
}_{n_2\;terms})\stackrel{w}\succeq
(\underbrace{\mu_1,\mu_1,...,\mu_1}_{n_1\;terms},\underbrace{\eta ,\eta
,...,\eta }_{n_2\;terms})\Rightarrow X_{1:n}\gtrsim_{hr}Y_{1:n}.$$
\end{thm}
\textbf{Proof:} Note that
\begin{eqnarray*}\frac{r_{X_{1:n}}(x)}{r_{Y_{1:n}}(x)}&=&\frac{
\frac{n_1}{1-\bar{\lambda_1}\bar{F}(x)}+\frac{n_2}{1-\bar{\eta}\bar{F}(x)}}{
\frac{n_1}{1-\bar{\mu_1}\bar{F}(x)}+\frac{n_2}{1-\bar{\eta}\bar{F}(x)}}\\
&=&\gamma(x),~
say.\end{eqnarray*}
We need to prove that $\gamma(x)$
is decreasing in $x>0$. As earlier, let us take
$u(x)=1/(1-x)$ and $v(x)=x/(1-x)$, which are increasing in $x$. Now
differentiating $\gamma(x)$ with respect to $x$, we have
\begin{eqnarray*}\gamma'(x)&\stackrel{sign}{=}&n_1^2
u(\bar{\lambda_1}\bar{F}(x))u(\bar{\mu_1}\bar{F}(x))[v(\bar{\mu_1}\bar{F}(x))-v(\bar{\lambda_1}\bar{F}(x))]+n_1
n_2
u(\bar{\lambda_1}\bar{F}(x))u(\bar{\eta}\bar{F}(x))\\&&[v(\bar{\eta}\bar{F}(x))-v(\bar{\lambda_1}\bar{F}(x))]+n_1
n_2
u(\bar{\eta}\bar{F}(x))u(\bar{\mu_1}\bar{F}(x))[v(\bar{\mu_1}\bar{F}(x))-v(\bar{\eta}\bar{F}(x))]\\&=&\psi(x),
say.\end{eqnarray*} Now the conditions $\lambda_1\leq
\min\{\eta,\mu_1\}$
 and $(\underbrace{\lambda_1,\lambda_1,...,\lambda_1}_{n_1\;terms},\underbrace{\eta
,\eta ,...,\eta }_{n_2\;terms})\stackrel{w}\succeq
(\underbrace{\mu_1,\mu_1,...,\mu_1}_{n_1\;terms},\underbrace{\eta
,\eta ,...,\eta }_{n_2\;terms})$ together is equivalent to the fact
that
$\lambda_1\leq \eta \leq \mu_1$ or $\lambda_1\leq \mu_1\leq \eta$.\\
Case I: Let $\lambda_1\leq \eta \leq \mu_1$. Then $\psi(x)\leq 0$.\\
Case II: Let $\lambda_1\leq \mu_1\leq \eta$. Then we have
$$u(\bar{\lambda_1}\bar{F}(x))\geq u(\bar{\mu_1}\bar{F}(x))\geq
u(\bar{\eta}\bar{F}(x))$$
 and
 $$v(\bar{\lambda_1}\bar{F}(x))\geq
v(\bar{\mu_1}\bar{F}(x))\geq v(\bar{\eta}\bar{F}(x)),$$ so that
\begin{eqnarray*}\psi(x)&\leq& n_1^2
u(\bar{\lambda_1}\bar{F}(x))u(\bar{\mu_1}\bar{F}(x))[v(\bar{\mu_1}\bar{F}(x))-v(\bar{\lambda_1}\bar{F}(x))]\\&&+n_1
n_2
u(\bar{\lambda_1}\bar{F}(x))u(\bar{\eta}\bar{F}(x))[v(\bar{\mu_1}\bar{F}(x))-v(\bar{\lambda_1}\bar{F}(x))]]\\&=&
n_1
u(\bar{\lambda_1}\bar{F}(x))[v(\bar{\mu_1}\bar{F}(x))-v(\bar{\lambda_1}\bar{F}(x))][n_1
u(\bar{\mu_1}\bar{F}(x))+n_2 u(\bar{\eta}\bar{F}(x))]\\&\leq& 0.
\end{eqnarray*} Hence the theorem follows.$\hfill\Box$
%Case III: Let $\eta\leq\lambda_1\leq \mu_1$. Then
%we have
%$$u(\bar{\eta}\bar{F}(x))\geq
%u(\bar{\lambda_1}\bar{F}(x))\geq u(\bar{\mu_1}\bar{F}(x))$$ and
%$$v(\bar{\eta}\bar{F}(x))\geq v(\bar{\lambda_1}\bar{F}(x))\geq
%v(\bar{\mu_1}\bar{F}(x)),$$ so that
%\begin{eqnarray*}\psi(x)&\leq& \textcolor[rgb]{0.98,0.00,0.00}{n_1^2
%u(\bar{\lambda_1}\bar{F}(x))u(\bar{\mu_1}\bar{F}(x))[v(\bar{\mu_1}\bar{F}(x))-v(\bar{\lambda_1}\bar{F}(x))]}\\&&\textcolor[rgb]{0.98,0.00,0.00}{+n_1
%n_2
%u(\bar{\lambda_1}\bar{F}(x))u(\bar{\eta}\bar{F}(x))[v(\bar{\mu_1}\bar{F}(x))-v(\bar{\lambda_1}\bar{F}(x))]]}\\&=&
%n_1
%u(\bar{\lambda_1}\bar{F}(x))[v(\bar{\mu_1}\bar{F}(x))-v(\bar{\lambda_1}\bar{F}(x))][n_1
%u(\bar{\mu_1}\bar{F}(x))+n_2 u(\bar{\eta}\bar{F}(x))]\\&\leq&
%0.\end{eqnarray*}

By taking $n_1=n_2=1$, the following corollary immediately follows
from Theorem \ref{thag3}.
\begin{coro}Let $X_1$ and $X_2$ be independent following the PO model with parameters $\lambda_1$ and $\eta$ respectively, and let $Y_1$ and $Y_2$ be independent following the PO model with parameters $\mu_1$ and $\eta$ respectively.
Then for $\lambda_1\leq \min\{\eta,\mu_1\}$,
$$(\lambda_1,\eta)\stackrel{w}\succeq(\mu_1,\eta)\Rightarrow X_{1:2}\gtrsim_{hr}Y_{1:2}.$$\end{coro}
The following lemma, required to prove the next theorem, has been
borrowed from \cite{kund}.
\begin{lem}\label{lema}If $\lambda_1\geq \mu_1\geq
\mu_2\geq \lambda_2$ or $\lambda_1\leq \mu_1\leq \mu_2\leq
\lambda_2$, and $n_1 \lambda_1+n_2 \lambda_2=n_1 \mu_1+n_2 \mu_2$,
then
$$(\underbrace{\lambda_1,\lambda_1,...,\lambda_1}_{n_1\;terms},\underbrace{\lambda_2,\lambda_2,...,\lambda_2}_{n_2\;terms})\stackrel{m}\succeq
(\underbrace{\mu_1,\mu_1,...,\mu_1}_{n_1\;terms},\underbrace{\mu_2,\mu_2,...,\mu_2}_{n_2\;terms}).$$
\end{lem}
The following theorem shows that under a different kind of
restriction on the model parameters than what is given in Theorem
\ref{thag3}, the condition of majorization order in Theorem
\ref{thag} can be replaced by the weak supermajorization order.
\begin{thm}\label{thag4}Let both $X$ and $Y$ follow the multiple-outlier PO
model with $X_i\sim PO(\bar{F},\lambda_1)$, $Y_i\sim
PO(\bar{F},\mu_1)$, for $i=1,2,...,n_1$, $X_j\sim
PO(\bar{F},\lambda_2)$, $Y_j\sim PO(\bar{F},\mu_2)$, for
$j=n_1+1,n_1+2,...,n_1+n_2(=n)$. Then, for \{$\lambda_1\leq
\mu_1\leq \mu_2\leq \lambda_2$\} or \{$\lambda_1\geq \mu_1\geq
\mu_2\geq \lambda_2$\},
$$(\underbrace{\lambda_1,\lambda_1,...,\lambda_1}_{n_1\;terms},\underbrace{\lambda_2,\lambda_2,...,\lambda_2}_{n_2\;terms})\stackrel{w}\succeq
(\underbrace{\mu_1,\mu_1,...,\mu_1}_{n_1\;terms},\underbrace{\mu_2,\mu_2,...,\mu_2}_{n_2\;terms})\Rightarrow
X_{1:n}\gtrsim_{hr} Y_{1:n}.$$
\end{thm}
\textbf{Proof:} Suppose that the first set of conditions holds. The
weak supermajorization order gives that $\lambda_1\leq \mu_1$ and
$n_1 \lambda_1+r \lambda_2\leq n_1 \mu_1+r \mu_2$, for
$r=1,2,...,n_2$. If $n_1 \lambda_1+n_2 \lambda_2=n_1 \mu_1+n_2
\mu_2$ holds then, under the given condition, the result follows
from Theorem \ref{thag}. Suppose that $n_1 \lambda_1+n_2
\lambda_2<n_1 \mu_1+n_2 \mu_2$. Then there exists an $\eta$
satisfying $\lambda_1< \eta \leq \mu_1$ such that $n_1 \eta+n_2
\lambda_2=n_1 \mu_1+n_2 \mu_2$. Let $X_{1:n}^{*}$ be the lifetime of
a series system formed by $n$ components having lifetimes
$X_1^{*},X_2^{*},...,X_n^{*}$, where $X_i^{*}\sim PO(\bar{F},\eta)$,
for $i=1,2,...,n_1$ and $X_j^{*}\sim PO(\bar{F},\lambda_2)$, for
$j=n_1+1,n_1+2,...,n_1+n_2(=n)$. Then, from Lemma \ref{lema} and
Theorem \ref{thag}, we have $X_{1:n}^{*}\gtrsim_{hr}Y_{1:n}$.
Further, we have $\lambda_1< \eta \leq \lambda_2$ and
$$(\underbrace{\lambda_1,\lambda_1,...,\lambda_1}_{n_1\;terms},\underbrace{\lambda_2,\lambda_2,...,\lambda_2}_{n_2\;terms})\stackrel{w}\succeq(\underbrace{\eta,\eta,...,\eta}_{n_1\;terms},\underbrace{\lambda_2,\lambda_2,...,\lambda_2}_{n_2\;terms}).$$
So, from Theorem \ref{thag3}, it follows that $X_{1:n}\gtrsim_{hr}
X_{1:n}^{*}$. Hence $X_{1:n}\gtrsim_{hr} Y_{1:n}$. The proof for the
second set of conditions can be done in a similar way.$\hfill\Box$

By taking $n_1=n_2=1$, the following corollary immediately follows
from Theorem \ref{thag4}.
\begin{coro}Let $X_1$ and $X_2$ be independent following PO model with parameters $\lambda_1$ and $\lambda_2$ respectively, and let $Y_1$ and $Y_2$ be independent following PO model with parameters $\mu_1$ and $\mu_2$ respectively.
Then
$$(\lambda_1,\lambda_2)\stackrel{w}\succeq(\mu_1,\mu_2)\Rightarrow X_{1:2}\gtrsim_{hr}Y_{1:2},$$
where $\{\lambda_1\leq \mu_1\leq \mu_2\leq\lambda_2\}$ or
$\{\lambda_1\geq \mu_1\geq \mu_2\geq
\lambda_2\}$.$\hfill\Box$\end{coro}
The following theorem shows
that, under certain condition, a series system with homogeneous
components ages faster than that with heterogeneous ones in terms of
the hazard rate.
\begin{thm}\label{thhma}Suppose lifetime vectors $X\sim PO(\bar{F},\boldsymbol\lambda)$ and $Y\sim
PO(\bar{F},\lambda \boldsymbol 1)$. Then,
$X_{1:n}\gtrsim_{hr}Y_{1:n}$ if $\lambda \geq
\frac{1}{n}\sum_{i=1}^n \lambda_i$.
\end{thm}
\textbf{Proof:} We have
$$\frac{r_{X_{1:n}}(x)}{r_{Y_{1:n}}(x)}=\frac{1-\bar{\lambda} \bar{F}(x)}{n}\sum_{i=1}^n \frac{1}{1-\bar{\lambda}_{i}\bar{F}(x)}.$$
Now, differentiating the above expression with respect to $x$, we
have, for $x>0$,
$$\frac{d}{dx}\left(\frac{r_{X_{1:n}}(x)}{r_{Y_{1:n}}(x)}\right)=\frac{f(x)(1-\bar{\lambda} \bar{F}(x))}{n}\left[\left(\frac{\bar{\lambda}}{1-\bar{\lambda} \bar{F}(x)}\right)\left(\sum_{i=1}^n \frac{1}{1-\bar{\lambda}_{i} \bar{F}(x)}\right)-\sum_{i=1}^n \frac{\bar{\lambda}_{i}}{(1-\bar{\lambda}_{i} \bar{F}(x))^2}\right],$$
so that $\frac{r_{X_{1:n}}(x)}{r_{Y_{1:n}}(x)}$ is decreasing if
\begin{equation}\label{aaa}
\left(\frac{\bar{\lambda} \bar{F}(x)}{1-\bar{\lambda}
\bar{F}(x)}\right)\left(\sum_{i=1}^n \frac{1}{1-\bar{\lambda}_{i}
\bar{F}(x)}\right)\leq \sum_{i=1}^n \frac{\bar{\lambda}_{i}
\bar{F}(x)}{(1-\bar{\lambda}_{i} \bar{F}(x))^2}.
\end{equation} From
Ceby\v{s}ev's inequality (cf. \citealp{mitr}, p. 240), (\ref{aaa}) holds if
\begin{equation*}\left(\frac{\bar{\lambda}
\bar{F}(x)}{1-\bar{\lambda} \bar{F}(x)}\right)\left(\sum_{i=1}^n
\frac{1}{1-\bar{\lambda}_{i} \bar{F}(x)}\right)\leq
\frac{1}{n}\left(\sum_{i=1}^n \frac{\bar{\lambda}_{i}
\bar{F}(x)}{1-\bar{\lambda}_{i} \bar{F}(x)}\right)\left(\sum_{i=1}^n
\frac{1}{1-\bar{\lambda}_{i} \bar{F}(x)}\right)\end{equation*} or
equivalently, \begin{equation}\label{ineq}\frac{\bar{\lambda}
\bar{F}(x)}{1-\bar{\lambda} \bar{F}(x)}\leq \frac{1}{n}\sum_{i=1}^n
\frac{\bar{\lambda}_{i} \bar{F}(x)}{1-\bar{\lambda}_{i}
\bar{F}(x)}.\end{equation} Let $\phi(x)=x/(1-x)$, which is
increasing and convex in $x$. Now (\ref{ineq}) holds if
$$\phi(\bar{\lambda}
\bar{F}(x))\leq \frac{1}{n}\sum_{i=1}^n \phi(\bar{\lambda}_{i}
\bar{F}(x)),$$ i.e. if
$$\phi(\bar{\lambda}
\bar{F}(x))\leq \phi\left(\frac{1}{n}\sum_{i=1}^n\bar{\lambda}_{i}
\bar{F}(x)\right),$$ which follows from the fact that $\phi$ is
convex. Now the theorem holds because $\phi$ is increasing.$\hfill\Box$\\

In case of multiple-outlier model, below we study the likelihood
ratio ordering between two series systems with heterogeneous
components. The result under majorization order follows from
Theorems \ref{thhr} and \ref{thag}, whereas the result under weak
supermajorization order follows from Theorems \ref{thhr} and
\ref{thag4}, by using the fact that the hazard rate order together
with the relative ageing order in the sense of hazard rate implies
the likelihood ratio order.
\begin{thm}\label{thlr}Let both $X$ and $Y$ follow the multiple-outlier PO
model such that $X_i\sim PO(\bar{F},\lambda_1)$, $Y_i\sim
PO(\bar{F},\mu_1)$, for $i=1,2,...,n_1$, $X_j\sim
PO(\bar{F},\lambda_2)$, $Y_j\sim PO(\bar{F},\mu_2)$, for
$j=n_1+1,n_1+2,...,n_1+n_2(=n)$. Then
$$(\underbrace{\lambda_1,\lambda_1,...,\lambda_1}_{n_1\;terms},\underbrace{\lambda_2,\lambda_2,...,\lambda_2}_{n_2\;terms})\stackrel{m}\succeq
(\text{resp.} \stackrel{w}\succeq)~
(\underbrace{\mu_1,\mu_1,...,\mu_1}_{n_1\;terms},\underbrace{\mu_2,\mu_2,...,\mu_2}_{n_2\;terms})\Rightarrow
X_{1:n}\leq_{lr}Y_{1:n},$$ provided $\{(\lambda_1,\lambda_2)\in
\mathcal{E_{+}},(\mu_1,\mu_2)\in \mathcal{E_{+}}\}$ or
$\{(\lambda_1,\lambda_2)\in \mathcal{D_{+}},(\mu_1,\mu_2)\in
\mathcal{D_{+}}\}\;($resp. \{$\lambda_1\leq \mu_1\leq \mu_2\leq
\lambda_2$\} or \{$\lambda_1\geq \mu_1\geq \mu_2\geq \lambda_2\})$
holds.$\hfill\Box$
\end{thm}
The following theorem gives a condition under which a series system
with homogeneous components and that with heterogeneous ones are
ordered in terms of the likelihood ratio order. The proof follows
from Theorem \ref{thhma} and Corollary \ref{corhr}.
\begin{thm}\label{thlr1}Suppose lifetime vectors $X\sim PO(\bar{F},\boldsymbol\lambda)$ and $Y\sim
PO(\bar{F},\lambda \boldsymbol 1)$. Then, $X_{1:n}\leq_{lr}Y_{1:n}$
if $\lambda \geq \frac{1}{n}\sum_{i=1}^n \lambda_i$.
\end{thm}
\section{Parallel systems with component lifetimes following the PO model}
In this section we compare lifetimes of two parallel systems of
heterogeneous components having lifetimes following the PO model
with respect to some stochastic orders. We also compare lifetimes of
two parallel systems, one comprising of heterogeneous components and
another of homogeneous components.\par %Let $X=(X_1,X_2,...,X_n)$ and
%$Y=(Y_1,Y_2,...,Y_n)$ be two sets of independent random variables,
%each following PO model. Let $X\sim PO(\bar{F},\boldsymbol\lambda)$
%and $Y\sim PO(\bar{F},\boldsymbol\mu)$, where $\bar{F}$ is the
%baseline survival function,
%$\boldsymbol\lambda=(\lambda_1,\lambda_2,...,\lambda_n)$ and
%$\boldsymbol\mu=(\mu_1,\mu_2,...,\mu_n)$, $\lambda_{i}>0$ and
%$\mu_i>0$, $i=1,2,...,n$.
We have the survival functions of $X_{n:n}$ and $Y_{n:n}$,
respectively, as
\begin{equation}\label{eqstp}\bar{F}_{X_{n:n}}(x)=1-\prod_{i=1}^{n}
(1-\bar{F}_{X_i}(x))=1-\prod_{i=1}^{n}\left(\frac{1-
\bar{F}(x)}{1-\bar{\lambda}_i\bar{F}(x)}\right),
\end{equation} and
$$\bar{F}_{Y_{n:n}}(x)=1-\prod_{i=1}^{n}
(1-\bar{F}_{Y_i}(x))=1-\prod_{i=1}^{n}\left(\frac{1-
\bar{F}(x)}{1-\bar{\mu}_i\bar{F}(x)}\right),$$ where
$\bar{\lambda}_i=1-\lambda_i$ and $\bar{\mu}_i=1-\mu_i$, for
$i=1,2,\ldots,n$. Also the reversed hazard rate functions of
$X_{n:n}$ and $Y_{n:n}$ are obtained, respectively, as
\begin{equation}\label{eqhrp}\tilde{r}_{X_{n:n}}(x)=\sum_{i=1}^{n}
\tilde{r}_{X_i}(x)=\sum_{i=1}^n \frac{\lambda_i
\tilde{r}(x)}{1-\bar{\lambda}_i\bar{F}(x)},
\end{equation} and
$$\tilde{r}_{Y_{n:n}}(x)=\sum_{i=1}^{n} \tilde{r}_{Y_i}(x)=\sum_{i=1}^n
\frac{\mu_i \tilde{r}(x)}{1-\bar{\mu}_i\bar{F}(x)}.$$ If $X\sim
PO(\bar{F},\lambda \boldsymbol 1)$, $\lambda>0$, then the survival
function and the reversed hazard rate function of $X_{n:n}$ are given,
respectively, by
$$\bar{F}_{X_{n:n}}(x)=1-\left(\frac{
1-\bar{F}(x)}{1-\bar{\lambda}\bar{F}(x)}\right)^n,$$ and
$$\tilde{r}_{X_{n:n}}(x)=\frac{n \lambda
\tilde{r}(x)}{1-\bar{\lambda}\bar{F}(x)},$$ where
$\bar{\lambda}=1-\lambda$.

The following theorem compares the lifetimes of two parallel systems
formed out of $n$ heterogeneous components following PO model in
terms of reversed hazard rate order.
\begin{thm}\label{thprh} Suppose that lifetime vectors $X\sim PO(\bar{F},\boldsymbol\lambda)$ and $Y\sim
PO(\bar{F},\boldsymbol\mu)$. Then
$$\boldsymbol\lambda\stackrel{w}\succeq\boldsymbol\mu~\text{implies}~ X_{n:n}\leq_{rhr}Y_{n:n}.$$ \end{thm}
\textbf{Proof:} Differentiating (\ref{eqhrp}) with respect to
$\lambda_i$ we have \begin{eqnarray*}\frac{\partial
\tilde{r}_{X_{n:n}}(x)}{\partial
\lambda_i}&=&\frac{\tilde{r}(x)(1-\bar{F}(x))}{(1-\bar{\lambda}_i\bar{F}(x))^
2}\\&\geq& 0,\end{eqnarray*} so that $\tilde{r}_{X_{n:n}}(x)$ is
increasing in each $\lambda_i$, $i=1,2,...,n.$ Also
$\tilde{r}_{X_{n:n}}(x)$ is symmetric with respect to
$(\lambda_1,\lambda_2,...,\lambda_n)\in \mathbb{R}^n$. For $1\leq
i\leq j\leq n$,
\begin{eqnarray*}(\lambda_i-\lambda_j)\left(\frac{\partial \tilde{r}_{X_{n:n}}(x)}{\partial
\lambda_i}-\frac{\partial \tilde{r}_{X_{n:n}}(x)}{\partial
\lambda_j}\right)&=&(\lambda_i-\lambda_j)\tilde{r}(x)(1-\bar{F}(x))\left[\frac{1}{(1-\bar{\lambda}_i\bar{F}(x))^
2}-\frac{1}{(1-\bar{\lambda_j}\bar{F}(x))^
2}\right]\\&\stackrel{sign}{=}&(\lambda_i-\lambda_j)\left((1-\bar{\lambda_j}\bar{F}(x))^2-(1-\bar{\lambda}_i\bar{F}(x))^
2\right)\\&\leq& 0.
\end{eqnarray*}
So, from Lemma \ref{le2a}, it follows that $\tilde{r}_{X_{n:n}}(x)$
is Schur-concave in
$\boldsymbol\lambda=(\lambda_1,\lambda_2,...,\lambda_n)\in
\mathbb{R}^n$. Thus, from Lemma \ref{le3}, we have
$\tilde{r}_{X_{n:n}}(x)\leq \tilde{r}_{Y_{n:n}}(x)$, for all $x$,  whenever
$\boldsymbol\lambda\stackrel{w}\succeq\boldsymbol\mu$.$\hfill\Box$

Since
$(\lambda_1,\lambda_2,\ldots,\lambda_n)\stackrel{w}\succeq(\underbrace{\lambda,\lambda,\ldots,\lambda}_{n\;terms})$,
for $\lambda\geqslant\frac{1}{n}\sum_{i=1}^n\lambda_i$, the
following corollary immediately follows from the above theorem.
\begin{coro}\label{corprh}Suppose lifetime vectors $X\sim PO(\bar{F},\boldsymbol\lambda)$ and $Y\sim
PO(\bar{F},\lambda \boldsymbol 1)$. Then, $X_{n:n}\leq_{rhr}Y_{n:n}$
if $\lambda \geq \frac{1}{n}\sum_{i=1}^n \lambda_i$.
$\hfill\Box$
\end{coro}
One may wonder whether the condition of weakly supermajorization
order can be replaced by $p$-larger order. This is answered in
negative in Counterexample \ref{nstupl}  where it is shown that,
even for usual stochastic order, the condition of weakly
supermojorization order given in the above theorem cannot be
replaced by $p$-larger order.$\hfill\Box$

If two parallel systems are formed -- one out of heterogeneous
components under the PO model and the other of homogeneous
components, then the condition under which the former dominates the
latter in usual stochastic order is discussed in the following
theorem.
\begin{thm}\label{thpstgm}Suppose that lifetime vectors $X\sim PO(\bar{F},\boldsymbol\lambda)$ and $Y\sim
PO(\bar{F},\lambda \boldsymbol 1)$. Then $X_{n:n}\geq_{st}Y_{n:n}$
if $\lambda=\sqrt[n]{\lambda_1 \lambda_2 \cdots
\lambda_n}$.\end{thm} \textbf{Proof:} Write
\begin{eqnarray*}\bar{F}_{X_{n:n}}(x)=
\phi(\lambda_1,\lambda_2,...,\lambda_n)\end{eqnarray*} Then we have
$$\frac{\partial\phi}{\partial
\lambda_i}=\bar{F}(x)[1-\bar{F}_{X_{n:n}}(x)]\frac{1}{1-\bar{\lambda}_i\bar{F}(x)}.$$
Let $\lambda_p=\min_{1\leq i\leq n}\lambda_i$ and
$\lambda_q=\max_{1\leq i\leq n}\lambda_i$. Then
\begin{eqnarray*}
\left(\frac{1}{\prod_{i\neq p}\lambda_i}\right)\frac{\partial\phi}{\partial\lambda_p}-\left(\frac{1}{\prod_{i\neq q}
\lambda_i}\right)\frac{\partial\phi}{\partial
\lambda_q}&\stackrel{sign}{=}&\left(\frac{1}{\prod_{i\neq
p}\lambda_i}\right)\frac{1}{1-\bar{\lambda_p}\bar{F}(x)}-\left(\frac{1}{\prod_{i\neq
q}\lambda_i}\right)\frac{1}{1-\bar{\lambda_q}\bar{F}(x)}\\&\stackrel{sign}{=}&\frac{\lambda_p}{1-\bar{\lambda_p}\bar{F}(x)}-\frac{\lambda_q}{1-\bar{\lambda_q}\bar{F}(x)}\\
&\stackrel{sign}{=}& (\lambda_p-\lambda_q)(1-\bar{F}(x))<0.
\end{eqnarray*}
So $(\lambda_p-\lambda_q)\left(\frac{1}{\prod_{i\neq p}
\lambda_i}\frac{\partial\phi}{\partial
\lambda_p}-\frac{1}{\prod_{i\neq q}
\lambda_i}\frac{\partial\phi}{\partial \lambda_q}\right)>0$. Thus,
from Lemma \ref{lemgm}, we have, for $\lambda=\sqrt[n]{\lambda_1
\lambda_2 \cdots \lambda_n}$,
$\phi(\lambda_1,\lambda_2,...,\lambda_n)\geq
\phi(\lambda,\lambda,...,\lambda)$, i.e.
$X_{n:n}\geq_{st}Y_{n:n}$.$\hfill\Box$\par

%Following counterexample shows that even
%under majorization order between $\boldsymbol\lambda$ and
%$\boldsymbol\mu$, $X_{n:n}$ and $Y_{n:n}$ may not be ordered in
%terms of likelihood ratio order.
%\begin{counterexample}\label{nlrp}
%Suppose lifetime vectors $X=(X_1,X_2,X_3)\sim
%PO(\bar{F},\boldsymbol\lambda)$ and $Y=(Y_1,Y_2,Y_3)\sim
%PO(\bar{F},\boldsymbol\mu)$, where the baseline survival function
%$\bar{F}(x)=e^{-2x}$, $\boldsymbol\lambda=(2,4.5,5)$ and
%$\boldsymbol\mu=(3,4,4.5)$. Here
%$\boldsymbol\lambda\succeq^{m}\boldsymbol\mu$. Now it is observed
%that $$f_{Y_{3:3}}(x)/f_{X_{3:3}}(x)=\frac{\sum_{i=1}^3
%\frac{\mu_i}{1-\bar{\mu}_i\bar{F}(x)}}{\sum_{i=1}^3
%\frac{\lambda_i}{1-\bar{\lambda_i}\bar{F}(x)}}
%\prod_{i=1}^3\frac{1-\bar{\lambda_i}\bar{F}(x)}{1-\bar{\mu}_i\bar{F}(x)}$$
%is nonmonotone.
%\end{counterexample}

If two parallel systems are formed out of heterogeneous components
satisfying the PO model then one may expect in the line of Theorem
\ref{thag} that there exists relative ageing in terms of reversed
hazard rate of the two systems whenever there is a majorization
order among the parameters of the two systems. However,
Counterexample \ref{nlrmop} shows that in case of multiple-outlier
model, under the majorization order, two parallel systems of
heterogeneous components may not be ordered with respect to relative
ageing in terms of reversed hazard rate.
\begin{remark}\label{bbb}Taking the random variables as in Counterexample \ref{nlrmop}, we see from Figure \ref{fig:notlrpo}(b) that
$f_{Y_{6:6}}(x)/f_{X_{6:6}}(x)$ is also non-monotone. This gives
that, in case of multiple-outlier model, under the majorization
order, two parallel systems with heterogeneous components may not be
ordered with respect to likelihood ratio order. \hfill$\Box$\par
\end{remark}
The Counterexample \ref{nlrpr} shows that a parallel system of
heterogeneous components may not be comparable with that of
homogeneous components with respect to relative ageing in terms of
reversed hazard rate, with or without the condition in Corollary
\ref{corprh}. That is, in case of parallel system, we cannot find
similar result in line of Theorem \ref{thhma}.

In case of multiple-outlier model, following theorem gives a
condition under which $X_{n:n}$ ages faster than $Y_{n:n}$ in terms
of reversed hazard rate.
\begin{thm}\label{thrhag}Let both $X$ and $Y$ follow the multiple-outlier PO
model such that $X_i\sim PO(\bar{F},\lambda_1)$, $Y_i\sim
PO(\bar{F},\mu_1)$, for $i=1,2,...,n_1$, $X_j\sim PO(\bar{F},\eta)$,
$Y_j\sim PO(\bar{F},\eta)$, for $j=n_1+1,n_1+2,...,n_1+n_2(=n)$.
Then
$$\lambda_1\leq \eta\leq \mu_1 \Rightarrow X_{n:n}\lesssim_{rhr}Y_{n:n}.$$
\end{thm}
\textbf{Proof:} Note that
\begin{equation*}\frac{\tilde{r}_{Y_{n:n}}(x)}{\tilde{r}_{X_{n:n}}(x)}=\frac{
\frac{n_1 \mu_1}{1-\bar{\mu_1}\bar{F}(x)}+\frac{n_2
\eta}{1-\bar{\eta}\bar{F}(x)}}{ \frac{n_1
\lambda_1}{1-\bar{\lambda_1}\bar{F}(x)}+\frac{n_2
\eta}{1-\bar{\eta}\bar{F}(x)}}=\gamma(x),~say.
\end{equation*}
We have to show that $\gamma(x)$ is
increasing in $x>0$. Let us write $u(x)=1/(1-x)$ and $v(x)=x/(1-x)$,
both of which are increasing in $x$. Now differentiating $\gamma(x)$
with respect to $x$, we have
\begin{eqnarray*}
    \gamma'(x)&\stackrel{sign}{=}&n_1^2
\lambda_1 \mu_1
u(\bar{\lambda_1}\bar{F}(x))u(\bar{\mu_1}\bar{F}(x)[v(\bar{\lambda_1}\bar{F}(x))-v(\bar{\mu_1}\bar{F}(x))]+n_1
n_2 \eta \mu_1
u(\bar{\mu_1}\bar{F}(x))u(\bar{\eta}\bar{F}(x))\\&&[v(\bar{\eta}\bar{F}(x))-v(\bar{\mu_1}\bar{F}(x))]+n_1
n_2\eta \lambda_1 u(\bar{\lambda_1}\bar{F}(x)
u(\bar{\eta}\bar{F}(x))[v(\bar{\lambda_1}\bar{F}(x))-v(\bar{\eta}\bar{F}(x))]\\&\geq&0,
\end{eqnarray*}
if $\lambda_1\leq \eta\leq \mu_1$. Hence the theorem follows.\hfill$\Box$

By taking $n_1=n_2=1$, the following corollary immediately follows
from the above theorem.
\begin{coro}
Let $X_1$ and $X_2$ be independent following the PO model with
parameters $\lambda_1$ and $\eta$ respectively, and let $Y_1$ and
$Y_2$ be independent following PO model with parameters $\mu_1$ and
$\eta$ respectively. Then
$$\lambda_1\leq
\eta\leq \mu_1 \Rightarrow X_{2:2}\lesssim_{rhr}Y_{2:2}.$$
\end{coro}
\begin{remark}
It is interesting to note that, the
    condition $\lambda_1\leqslant\eta\leqslant \mu_1$ is crucial for the
    above theorem to hold. The Counterexample \ref{nrhrag} shows that
    Theorem \ref{thrhag} does not hold if the given condition is
    replaced by $\lambda_1\leq \mu_1\leq \eta$. \hfill$\Box$
\end{remark}
If two parallel systems are formed out of components following the
multiple-outlier PO model, then one might be interested to know some
condition(s) under which these two models are comparable in terms of
likelihood ratio order. This is given in the next theorem, which may
be compared with Remark \ref{bbb}.
\begin{thm}\label{thpmolr}Let both $X$ and $Y$ follow the multiple-outlier PO
model such that $X_i\sim PO(\bar{F},\lambda_1)$, $Y_i\sim
PO(\bar{F},\mu_1)$, for $i=1,2,...,n_1$, $X_j\sim PO(\bar{F},\eta)$,
$Y_j\sim PO(\bar{F},\eta)$, for $j=n_1+1,n_1+2,...,n_1+n_2(=n)$.
Then
$$\lambda_1\leq \eta\leq \mu_1 \Rightarrow X_{n:n}\leq_{lr}Y_{n:n}.$$
\end{thm}
\textbf{Proof:} We need to show that
\begin{equation}
\frac{f_{Y_{n:n}}(x)}{f_{X_{n:n}}(x)}=\left(\frac{F_{Y_{n:n}}(x)}{F_{X_{n:n}}(x)}\right)\left( \frac{\tilde{r}_{Y_{n:n}}(x)}{\tilde{r}_{X_{n:n}}(x)}\right)\;{\rm is\; increasing\; in}\;x>0.
\end{equation}
We have $\lambda_1\leq \eta\leq \mu_1$, which
implies that
$(\underbrace{\lambda_1,\lambda_1,...,\lambda_1}_{n_1\;terms},\underbrace{\eta
,\eta ,...,\eta }_{n_2\;terms})\stackrel{w}\succeq
(\underbrace{\mu_1,\mu_1,...,\mu_1}_{n_1\;terms},\underbrace{\eta ,\eta
,...,\eta }_{n_2\;terms})$. So, from Theorem \ref{thprh}, under the given
condition, we get that $F_{Y_{n:n}}(x)/F_{X_{n:n}}(x)$ is increasing in $x>0$.
Again, Theorem \ref{thrhag} gives that, under the given condition,
$\tilde{r}_{Y_{n:n}}(x)/\tilde{r}_{X_{n:n}}(x)$ is increasing in
$x>0$. Hence the theorem follows.\hfill$\Box$

For $n_1=n_2=1$, the above theorem reduces to the following
corollary.
\begin{coro}Let $X_1$ and $X_2$ be independent following the PO model with parameters $\lambda_1$ and $\eta$ respectively, and let $Y_1$ and $Y_2$ be independent following PO model with parameters $\mu_1$ and $\eta$ respectively.
Then
$$\lambda_1\leq
\eta\leq \mu_1 \Rightarrow X_{2:2}\leq_{lr}Y_{2:2}.$$\end{coro}
\section{Examples and counterexamples}
In this section, we give some examples to illustrate the proposed
results, and some counterexamples are given wherever needed.
\subsection{Examples}
    In this subsection we give some examples to demonstrate the proposed results of this paper. The first example gives an application of Theorem~\ref{thst}.
    \begin{example}\label{exth3.1}
        Consider two series systems, each
        comprising of three components having lifetimes following the PO
        model with the common baseline survival function given by
        $\bar{F}(x)=e^{-(x/\beta)^k}$ with $\beta=0.4$, $k=2$, $x\geq 0$.
        Then the survival functions of two series systems are given by
        \begin{equation}\label{eqth3.1}\bar{F}_{X_{1:3}}(x)=\prod_{i=1}^{3}\frac{\lambda_i
            e^{-(x/0.4)^2}}{1-\bar{\lambda}_ie^{-(x/0.4)^2}}\qquad\text{and}\qquad
        \bar{F}_{Y_{1:3}}(x)=\prod_{i=1}^{3}\frac{\mu_i
            e^{-(x/0.4)^2}}{1-\bar{\mu}_ie^{-(x/0.4)^2}},\end{equation}
        respectively, where $(\lambda_1,\lambda_2,\lambda_3)=(2,3,5)$ and
        $(\mu_1,\mu_2,\mu_3)=(2.5,3.5,6)$ so that
        $(2,3,5)\stackrel{p}\succeq(2.5,3.5,6)$. In order to change the scale, we
        substitute $x=t/(1-t)$ in (\ref{eqth3.1}) so that, for
        $x\in[0,\infty)$, we have $t\in[0,1)$, and after this substitution,
        let us denote the expressions in (\ref{eqth3.1}) as $\xi_1(t)$ and
        $\xi_2(t)$, respectively. From Figure \ref{picth3.1} we observe
        that $\xi_1(t)\leq \xi_2(t)$ for all $t\in[0,1)$, which implies that
        $\bar{F}_{X_{1:3}}(x)\leq \bar{F}_{Y_{1:3}}(x)$ for all $x\geq 0$. Thus $X_{1:3}\leq_{st}Y_{1:3}.$ $\hfill\Box$
        \begin{figure}\begin{center}
                \includegraphics[width=12cm]{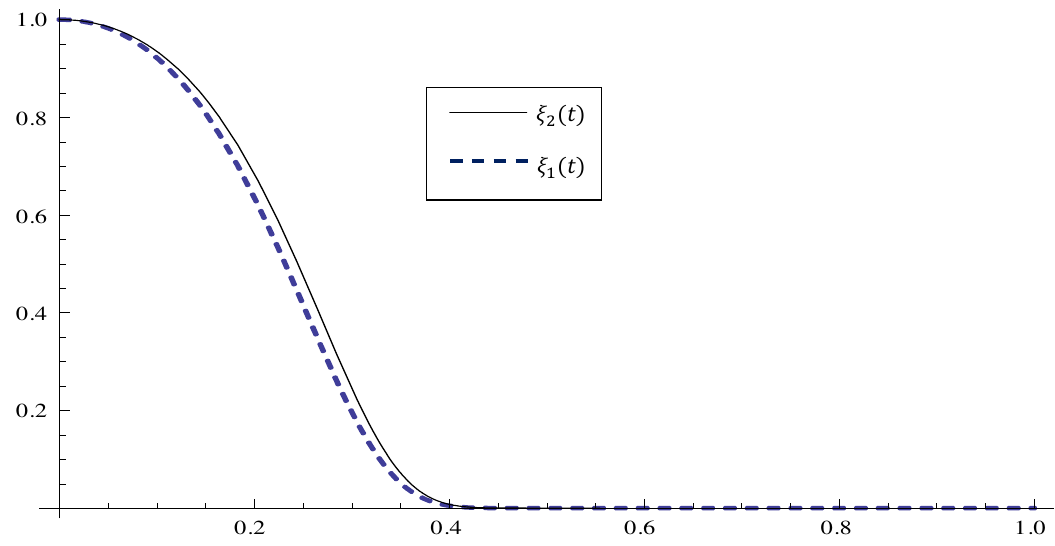}\\
                \caption{Plot of $\xi_1(t)$ and
                    $\xi_2(t)$ against $t\in[0,1]$.}\label{picth3.1}\end{center}
        \end{figure}
    \end{example}
    \hspace*{0.2 in}In the following example we illustrate the result given in Theorem~\ref{thhr}.
    \begin{example}\label{exth3.2}
        Consider two series systems, each
        comprising of three components having lifetimes following the PO
        model with the common baseline survival function given by
        $\bar{F}(x)=e^{-2x}$, $x\geq 0$. Then the survival functions of the two series systems
        are given by
        $$\bar{F}_{X_{1:3}}(x)=\prod_{i=1}^{3}\frac{\lambda_i
            e^{-2x}}{1-\bar{\lambda}_ie^{-2x}}\qquad\text{and}\qquad
        \bar{F}_{Y_{1:3}}(x)=\prod_{i=1}^{3}\frac{\mu_i
            e^{-2x}}{1-\bar{\mu}_ie^{-2x}},$$ respectively. Taking
        $(\lambda_1,\lambda_2,\lambda_3)=(3,4.5,6)$ and
        $(\mu_1,\mu_2,\mu_3)=(4,5,6)$ we observe that
        $(3,4.5,6)\stackrel{w}\succeq(4,5,6)$. Note that
        $$\frac{\bar{F}_{Y_{1:3}}(x)}{\bar{F}_{X_{1:3}}(x)}=\left(\frac{120}{81}\right)\left(\frac{1+5.5e^{-2x}+7e^{-4x}}{1+7e^{-2x}+12e^{-4x}}\right),$$
        which is increasing in $x\geq 0$, and hence $X_{1:3}\leq_{hr}Y_{1:3}.$ $\hfill\Box$
    \end{example}
    \hspace*{0.2 in}In the following example we demonstrate the result given in Theorem~\ref{thag}
    \begin{example}\label{exth3.3}
        Consider two series systems, each
        comprising of four components having lifetimes following the
        multiple-outlier PO model with the common baseline survival function given by
        $\bar{F}(x)=e^{-\frac{x^2}{2}}$, $x\geq 0$.
        Then the hazard rate functions of the two series systems are given by
        \begin{equation}
        \label{eqth3.3}r_{X_{1:4}}(x)= \frac{2r(x)}{1-\bar{\lambda_1}\bar{F}(x)}+
        \frac{2r(x)}{1-\bar{\lambda_2}\bar{F}(x)},
        \end{equation}
        and
        \begin{equation}\label{123}
        r_{Y_{1:4}}(x)=\frac{2r(x)}{1-\bar{\mu_1}\bar{F}(x)}+\frac{2r(x)}{1-\bar{\mu_2}\bar{F}(x)},
        \end{equation}
        respectively, where $r(\cdot)$ is the common hazard rate function of
        each of the components, and  $\lambda_1=2,\lambda_2=4.5,\mu_1=3,\mu_2=3.5$,
        so that $(2,2,4.5,4.5)\stackrel{m}\succeq(3,3,3.5,3.5)$. In order to
        change the scale, we substitute $x=t/(1-t)$ in (\ref{eqth3.3}) and (\ref{123}) so
        that, for $x\in[0,\infty)$, we have $t\in[0,1)$, and after this
        substitution, let us denote the expressions in (\ref{eqth3.3}) and (\ref{123}) as
        $l_1(t)$ and $l_2(t)$, respectively. From Figure~\ref{picth3.34}(a) we observe that $l_1(t)/l_2(t)$ is decreasing in
        $t\in[0,1)$, which is equivalent to the fact that $r_{X_{1:4}}(x)/r_{Y_{1:4}}(x)$ is decreasing in
        $x\geq 0$. Hence $X_{1:4}\gtrsim_{hr}Y_{1:4}.$ $\hfill\Box$
        \begin{figure}\begin{center}
                \includegraphics[width=12cm]{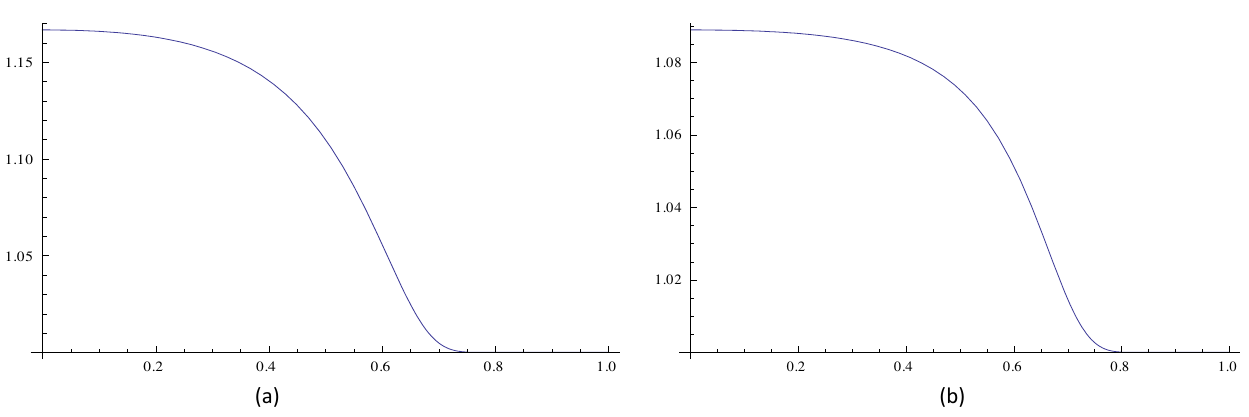}\\
                \caption{(a) Plot of $l_1(t)/l_2(t)$ against $t\in[0,1]$, (b) Plot of $\kappa_1(t)/\kappa_2(t)$ against $t\in[0,1]$.}\label{picth3.34}\end{center}
        \end{figure}
    \end{example}
    \hspace*{0.2 in}An application of Theorem~\ref{thag3} is given below.
    \begin{example}\label{exth3.5}
        Consider Example~\ref{exth3.3} with $\lambda_1=3,\lambda_2=\mu_2=\eta=4,\mu_1=3.5$,
        so that $(2,2,4,4)\stackrel{w}\succeq(3.5,3.5,4,4)$. After
        substituting $x=t/(1-t)$ in (\ref{eqth3.3}) and (\ref{123}), let us denote the
        expressions in (\ref{eqth3.3}) and (\ref{123}) as $\kappa_1(t)$ and $\kappa_2(t)$,
        respectively. From Figure \ref{picth3.34}(b) we observe that
        $\kappa_1(t)/\kappa_2(t)$ is decreasing in $t\in[0,1)$, which implies that
        $r_{X_{1:4}}(x)/r_{Y_{1:4}}(x)$ is decreasing in $x\geq 0$. Hence $X_{1:4}\gtrsim_{hr}Y_{1:4}.$ $\hfill\Box$
    \end{example}
    \hspace*{0.2 in}In the following example we demonstrate the result given in Theorem~\ref{thag4}.
    \begin{example}\label{exth3.6}
        Consider Example~\ref{exth3.3} with $\lambda_1=2,\lambda_2=4,\mu_1=3,\mu_2=3.5$, so
        that $(2,2,4,4)\stackrel{w}\succeq(3,3,3.5,3.5)$ and
        $\lambda_1<\mu_1<\mu_2<\lambda_2$. After substituting $x=t/(1-t)$ in
        (\ref{eqth3.3}) and (\ref{123}), let us denote the expressions in (\ref{eqth3.3}) and (\ref{123}) as
        $\zeta_1(t)$ and $\zeta_2(t)$, respectively. From Figure~\ref{picth3.6} we observe that $\zeta_1(t)/\zeta_2(t)$ is decreasing in
        $t\in[0,1)$, which is equivalent to the fact that $r_{X_{1:4}}(x)/r_{Y_{1:4}}(x)$ is decreasing in
        $x\geq 0$. Hence $X_{1:4}\gtrsim_{hr}Y_{1:4}.$ $\hfill\Box$
        \begin{figure}\begin{center}
                \includegraphics[width=12cm]{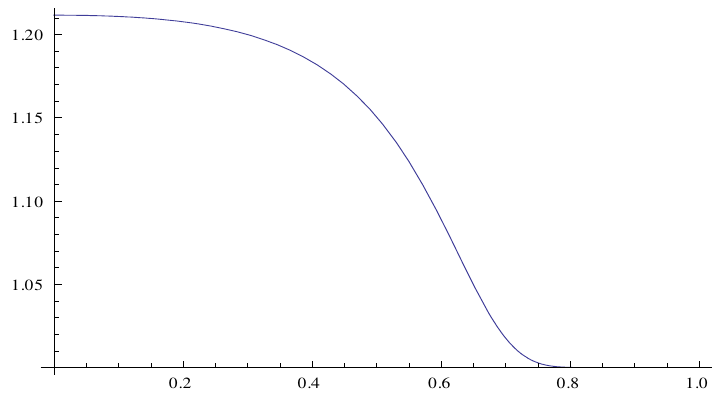}\\
                \caption{Plot of $\zeta_1(t)/\zeta_2(t)$ against $t\in[0,1]$.}\label{picth3.6}\end{center}
        \end{figure}
    \end{example}
    \hspace*{0.2 in}Below we give an example to illustrate the result given in Theorem~\ref{thprh}.
    \begin{example}\label{exth4.1}
        Consider two parallel systems, each
        comprising of three components having lifetimes following the PO
        model with the common baseline survival function given by
        $\bar{F}(x)=e^{-1.5x}$, $x\geq 0$. Then the reversed hazard rate functions of
        two parallel systems are given by
        \begin{equation}\label{eqth4.1}\tilde{r}_{X_{3:3}}(x)=\frac{1.5e^{-1.5x}}{1-e^{-1.5x}}\sum_{i=1}^3
        \frac{\lambda_i}{1-\bar{\lambda}_ie^{-1.5x}}
            \end{equation}
        and
        \begin{equation}\label{xyz}
        \tilde{r}_{Y_{3:3}}(x)=\frac{1.5e^{-1.5x}}{1-e^{-1.5x}}\sum_{i=1}^3
        \frac{\mu_i}{1-\bar{\mu}_ie^{-1.5x}},
        \end{equation} respectively,
        where $(\lambda_1,\lambda_2,\lambda_3)=(0.5,2.5,4)$ and
        $(\mu_1,\mu_2,\mu_3)=(1,3,5)$ so that
        $(0.5,2.5,4)\stackrel{w}\succeq(1,3,5)$. In order to change the scale, we
        substitute $x=t/(1-t)$ in (\ref{eqth4.1}) and (\ref{xyz}) so that, for
        $x\in[0,\infty)$, we have $t\in[0,1)$, and after this substitution,
        let us denote the expressions in (\ref{eqth4.1}) and (\ref{xyz}) as $\gamma_1(t)$
        and $\gamma_2(t)$, respectively. From Figure \ref{picth4.1} we
        observe that $\gamma_1(t)\leq \gamma_2(t),\;t\in[0,1)$, which
        implies that $\tilde{r}_{X_{3:3}}(x)\leq \tilde{r}_{Y_{3:3}}(x)$ for
        all $x\geq 0$. Hence $X_{3:3}\leq_{rhr}Y_{3:3}$. $\hfill\Box$
        \begin{figure}\begin{center}
                \includegraphics[width=12cm]{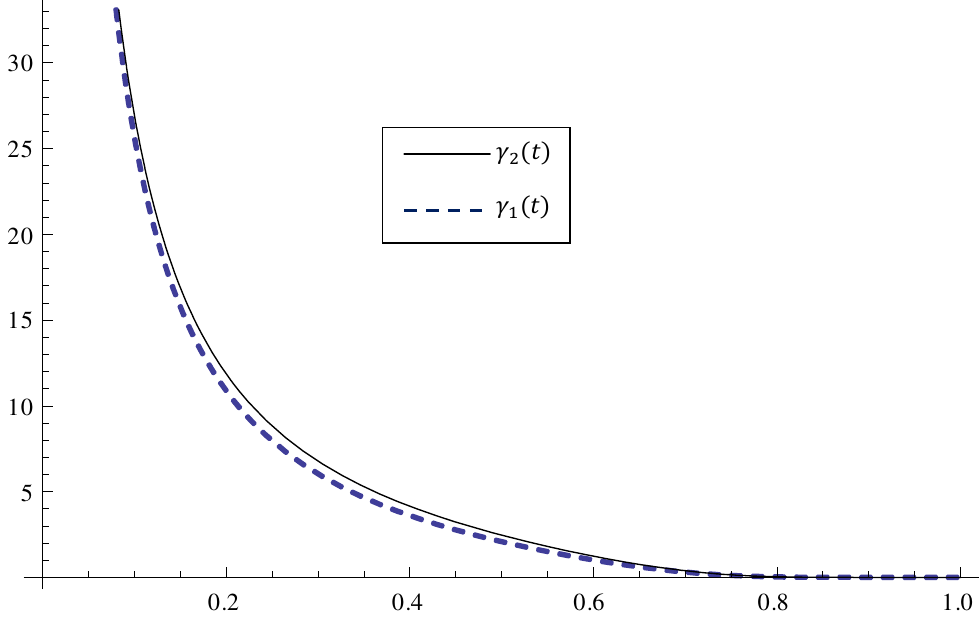}\\
                \caption{Plot of $\gamma_1(t)$ and
                    $\gamma_2(t)$ against $t\in[0,1]$.}\label{picth4.1}\end{center}
        \end{figure}
    \end{example}
    \hspace*{0.2 in}The following example demonstrates the results given in Theorems~\ref{thrhag}~and~\ref{thpmolr}.
    \begin{example}\label{exth4.3}
        Consider two parallel systems,
        each comprising of four components having lifetimes following the
        multiple-outlier PO model with the common baseline survival function given by
        $\bar{F}(x)=e^{-2x}$, $x\geq 0$. Then the  ratio of the reversed hazard
        rate functions of two parallel systems is given by
        \begin{equation}\frac{\tilde{r}_{Y_{4:4}}(x)}{\tilde{r}_{X_{4:4}}(x)}=\frac{
            \frac{2 \mu_1}{1-\bar{\mu_1}e^{-2x}}+\frac{2
                \eta}{1-\bar{\eta}e^{-2x}}}{ \frac{2
                \lambda_1}{1-\bar{\lambda_1}e^{-2x}}+\frac{2
                \eta}{1-\bar{\eta}e^{-2x}}},
        \end{equation} where
        $\lambda_1=2,\eta=3,\mu_1=4,$ so that $\lambda_1<\eta<\mu_1$. Note that $\tilde{r}_{Y_{4:4}}(x)/\tilde{r}_{X_{4:4}}(x)$ is increasing
        in $x\geq 0$. Hence $X_{4:4}\lesssim_{rhr}Y_{4:4}.$ Again, we have
        \begin{equation*}
        \frac{f_{Y_{4:4}}(x)}{f_{X_{4:4}}(x)}=\left(\frac{F_{Y_{4:4}}(x)}{F_{X_{4:4}}(x)}\right)\left(
        \frac{\tilde{r}_{Y_{4:4}}(x)}{\tilde{r}_{X_{4:4}}(x)}\right).\end{equation*}
        It can be verified that
        $F_{Y_{4:4}}(x)/F_{X_{4:4}}(x)=(1+e^{-2x})^2/(1+3e^{-2x})^2$ is increasing in $x\geq 0$, which
        implies that $f_{Y_{4:4}}(x)/f_{X_{4:4}}(x)$ is increasing in $x\geq 0$. Hence $X_{4:4}\leq_{lr}Y_{4:4}.$
\end{example}
\subsection{Counterexamples}\label{sp}
\begin{figure}[h]
    \begin{center}
  \includegraphics[width=12cm]{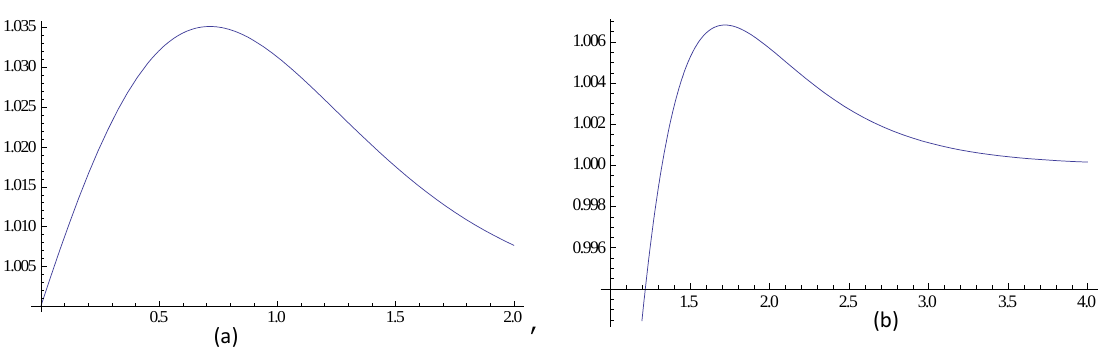}\\
  \caption{(a) Curve of $\tilde{r}_{Y_{6:6}}(x)/\tilde{r}_{X_{6:6}}(x)$   (b) Curve of $f_{Y_{6:6}}(x)/f_{X_{6:6}}(x)$} \label{fig:notlrpo}
\end{center}
\end{figure}

\begin{figure}[h]
    \begin{center}
        \includegraphics[width=13cm]{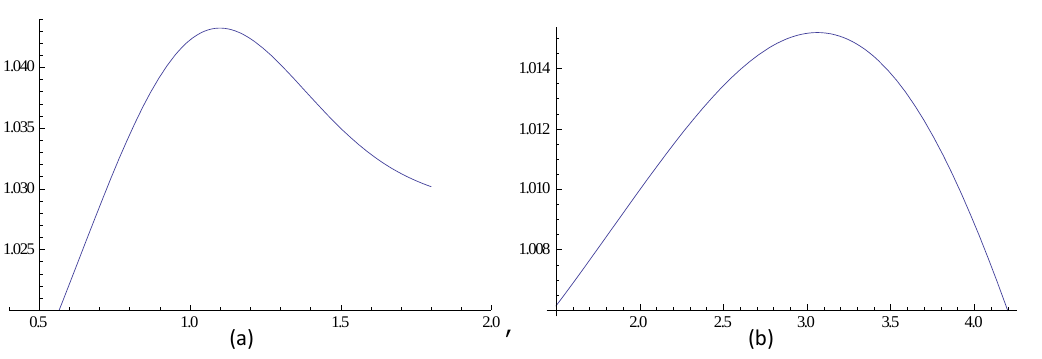}\\
        \caption{(a) Curve of $\tilde{r}_{Y_{4:4}}(x)/\tilde{r}_{X_{4:4}}(x)$ for $(\lambda_1,\lambda_2,\lambda_3,\lambda_4,\lambda,\beta,k)=(2,3,4,5,3.6,0.8,2)$, and
            (b) Curve of $\tilde{r}_{Y_{4:4}}(x)/\tilde{r}_{X_{4:4}}(x)$ for $(\lambda_1,\lambda_2,\lambda_3,\lambda_4,\lambda,\beta,k)=(2,3,4,5,3.4,3,2)$} \label{fig:notagpr}
        \end{center}
\end{figure}

A list of counterexamples are discussed in this subsection. The following counterexample shows that the $p$-larger order in Theorem \ref{thst} cannot be replaced by reciprocal majorization order.
\begin{counterexample}\label{nsturm}
Let $X=(X_1,X_2,X_3)$ and $Y=(Y_1,Y_2,Y_3)$ with $X_i\sim
PO(\bar{F},\lambda_i)$ and $Y_i\sim PO(\bar{F},\mu_i)$, $i=1,2,3$,
where the baseline survival function $\bar F$ is given by
$\bar{F}(x)=e^{-2 x},\;x>0$. Take
$(\lambda_1,\lambda_2,\lambda_3)=(2.2,3,5)$ and
$(\mu_1,\mu_2,\mu_3)=(2.8,3.2,3.3)$ so that
$(\lambda_1,\lambda_2,\lambda_3)\stackrel{rm}\succeq(\mu_1,\mu_2,\mu_3)$
but
$(\lambda_1,\lambda_2,\lambda_3)\stackrel{p}\nsucceq(\mu_1,\mu_2,\mu_3)$.
It is observed that, for $x=0.2$, $\bar{F}_{X_{1:3}}(x)=0.63929$ and
$\bar{F}_{Y_{1:3}}(x)=0.641646$. Again, for $x=0.8$,
$\bar{F}_{X_{1:3}}(x)=0.0861549$ and
$\bar{F}_{Y_{1:3}}(x)=0.084394$. So
$X_{1:3}\nleq_{st}Y_{1:3}$.\hfill$\Box$
\end{counterexample}
Below we show that weak majorization order in Theorem \ref{thhr} cannot be replaced by $p$-larger order.
\begin{counterexample}\label{nhrup}
Let $X=(X_1,X_2,X_3)$ and $Y=(Y_1,Y_2,Y_3)$ with $X_i\sim
PO(\bar{F},\lambda_i)$ and $Y_i\sim PO(\bar{F},\mu_i)$, $i=1,2,3$,
where the baseline survival function $\bar F$ is given by
$\bar{F}(x)=e^{-1.2 x},\;x>0$. Take
$(\lambda_1,\lambda_2,\lambda_3)=(2,3,5)$ and
$(\mu_1,\mu_2,\mu_3)=(2.8,3.2,3.4)$ so that
$(\lambda_1,\lambda_2,\lambda_3)\stackrel{p}\succeq(\mu_1,\mu_2,\mu_3)$
but
$(\lambda_1,\lambda_2,\lambda_3)\stackrel{w}\nsucceq(\mu_1,\mu_2,\mu_3)$.
It is observed that, for $x=0.2$, $r_{X_{1:3}}(x)=1.4273$ and
$r_{Y_{1:3}}(x)=1.3516$, and, for $x=1.8$, $r_{X_{1:3}}(x)=2.8722$
and $r_{Y_{1:3}}(x)=2.8907$. Thus, we have
$X_{1:3}\nleq_{hr}Y_{1:3}$.\hfill$\Box$
\end{counterexample}
In Theorem \ref{thag}, we have seen that, in the case of
multiple-outlier model, out of two series systems formed from
heterogeneous components, one may dominate the other in relative
ageing in terms of hazard rate, provided  the two sets of the
parameters of the model have majorization order between them.
However, this kind of result may not hold for parallel systems as we
see in the following counterexample.
\begin{counterexample}\label{nlrmop}
    Let $X=(X_1,X_2,...,X_6)$ and $Y=(Y_1,Y_2,...,Y_6)$, each follows the multiple-outlier PO
    model such that $X_i\sim PO(\bar{F},2)$, $Y_i\sim PO(\bar{F},3)$,
    for $i=1,2$, $X_j\sim PO(\bar{F},6)$, $Y_j\sim PO(\bar{F},5.5)$, for
    $j=3,4,5,6$, where the baseline survival function is given by
    $\bar{F}(x)=e^{-2x},\;x>0$. Clearly,
    $(2,2,6,6,6,6)\stackrel{m}\succeq (3,3,5.5,5.5,5.5,5.5)$. However,
    it is observed from Figure \ref{fig:notlrpo}(a) that
    $\tilde{r}_{Y_{6:6}}(x)/\tilde{r}_{X_{6:6}}(x)$ is non-monotone.\hfill$\Box$
\end{counterexample}
    That weak supermajorization order in Theorem \ref{thprh} cannot be replaced by $p$-larger order is shown in the following counterexample.
\begin{counterexample}\label{nstupl}
Let $X=(X_1,X_2,X_3)$ and $Y=(Y_1,Y_2,Y_3)$ with $X_i\sim
PO(\bar{F},\lambda_i)$ and $Y_i\sim PO(\bar{F},\mu_i)$, $i=1,2,3$,
where the baseline survival function is given by $\bar{F}(x)=e^{-1.8
x}$, $x>0$. Take $(\lambda_1,\lambda_2,\lambda_3)=(2,3,5)$ and
$(\mu_1,\mu_2,\mu_3)=(2.6,3.2,3.7)$ so that
$(\lambda_1,\lambda_2,\lambda_3)\stackrel{p}\succeq(\mu_1,\mu_2,\mu_3)$
but
$(\lambda_1,\lambda_2,\lambda_3)\stackrel{w}\nsucceq(\mu_1,\mu_2,\mu_3)$.
It is observed that, for $x=1.5$, $\bar{F}_{X_{3:3}}(x)=0.471629$
and $\bar{F}_{Y_{3:3}}(x)=0.459619$ so that
$X_{3:3}\nleq_{st}Y_{3:3}$. Now, take
$(\lambda_1,\lambda_2,\lambda_3)=(2.5,3,5)$ and
$(\mu_1,\mu_2,\mu_3)=(3,3.8,4.4)$ which give
$(\lambda_1,\lambda_2,\lambda_3)\stackrel{p}\succeq(\mu_1,\mu_2,\mu_3)$.
It is observed that, for $x=1.2$, $\bar{F}_{X_{3:3}}(x)=0.67176$ and
$\bar{F}_{Y_{3:3}}(x)=0.69449$ so that
$X_{3:3}\ngeq_{st}Y_{3:3}$.\hfill$\Box$
\end{counterexample}
%The following counterexample shows that the condition $\lambda\geqslant \sum_{i=1}^n\lambda_i/n$ in Corollary \ref{corprh} cannot be dropped.
\begin{counterexample}\label{nlrpr} Let $X=(X_1,X_2,X_3,X_4)$ and
$Y=(Y_1,Y_2,Y_3,Y_4)$, each follows the multiple-outlier PO model
such that $X_i\sim PO(\bar{F},\lambda_i)$, $i=1,2,3,4$ and $Y_i\sim
PO(\bar{F},\lambda)$, $i=1,2,3,4$, where the baseline survival
function is given by $\bar{F}(x)=e^{-(x/\beta)^k}$, $\beta,
k>0,\;x>0$. It is observed from Figure \ref{fig:notagpr}(a) that,
for $\lambda_1=2$, $\lambda_2=3$, $\lambda_3=4$, $\lambda_4=5$,
$\lambda=3.6$, $\beta=0.8$, and $k=2$,
$\tilde{r}_{Y_{4:4}}(x)/\tilde{r}_{X_{4:4}}(x)$ is non-monotone.
Again, for $\lambda_1=2$, $\lambda_2=3$, $\lambda_3=4$,
$\lambda_4=5$, $\lambda=3.4$, $\beta=3$ and $k=2$,
$\tilde{r}_{Y_{4:4}}(x)/\tilde{r}_{X_{4:4}}(x)$ is also non-monotone
as can be seen from Figure \ref{fig:notagpr}(b).\hfill$\Box$
\end{counterexample}
That the condition $\lambda_1\leqslant\eta\leqslant \mu_1$ in Theorem \ref{thrhag} cannot be dropped is shown below.
\begin{counterexample}\label{nrhrag}
Let $X_1$ and $X_2$ follow the PO model with parameters $\lambda_1$
and $\eta$ respectively, and let $Y_1$ and $Y_2$ follow the PO model
with parameters $\mu_1$ and $\eta$ respectively, where the baseline
distribution is exponential with parameter $\lambda=2$. Now, for
$\lambda_1=0.2$, $\mu_1=0.4$ and $\eta=0.9$,
$\tilde{r}_{Y_{2:2}}(x)/\tilde{r}_{X_{2:2}}(x)$ is non-monotone, as
we see from Figure \ref{fig:notrhrPr}.
\begin{figure}\begin{center}
  \includegraphics[width=12cm]{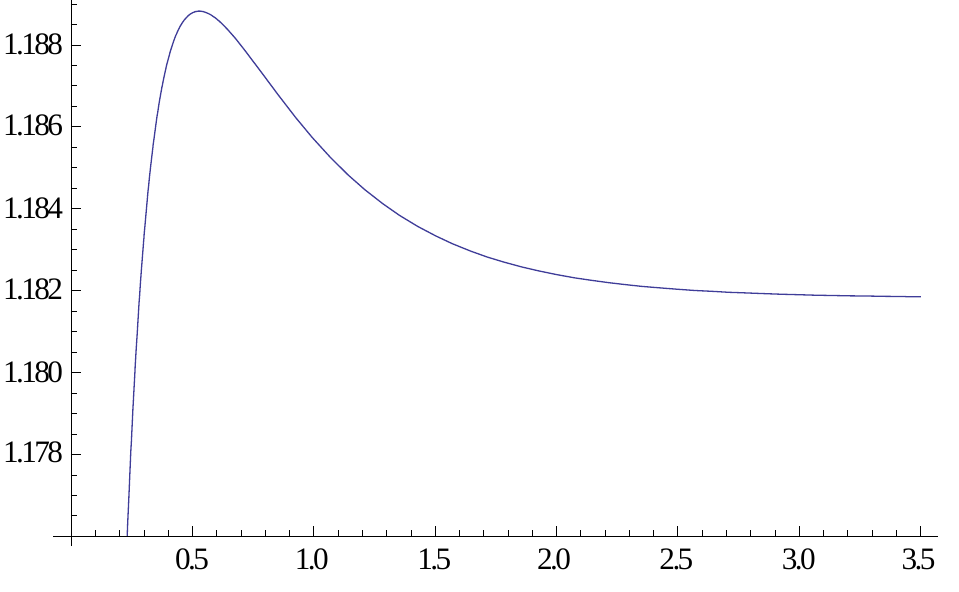}\\
  \caption{Curve of $\tilde{r}_{Y_{2:2}}(x)/\tilde{r}_{X_{2:2}}(x)$} \label{fig:notrhrPr}
\end{center}
\end{figure}
\end{counterexample}
\section{Conclusion}
In this paper, we have studied stochastic comparison of series and
parallel systems formed from independent heterogeneous components
having lifetimes following the PO model. Most of the results are
obtained using different concepts of majorization. We have also
compared a system formed of heterogeneous components with another
system of homogeneous components. We have derived conditions under
which two series systems with heterogeneous components are ordered
with respect to different stochastic orders; in the case of
multiple-outlier model, they are compared with respect to likelihood
ratio order and relative ageing in terms of hazard rate. We have
also derived conditions under which a series system with
heterogeneous components and that with homogeneous components are
ordered with respect to the above mentioned stochastic orderings. In
the case of parallel system, we have obtained conditions under which
two parallel systems with heterogeneous components are ordered with
respect to the usual stochastic order and reversed hazard rate
order. The comparison is also made in the case of a parallel system
with heterogeneous components and that with homogeneous components.
However, unlike series system, with suitable counterexamples, we
have shown that, even in the case of multiple-outlier model, under
majorization order, two parallel systems with heterogeneous
components may not be comparable with respect to likelihood ratio
order and relative ageing in terms of reversed hazard rate,
although, under more restricted conditions, we are able to compare
the parallel systems with respect to those stochastic orderings.

Similar kinds of results can be studied for a $k$-out-of-$n$ system
or equivalently, for $r$th largest order statistic (see the
explanation given in the introduction in this regard). It can be
noted that the expressions for different reliability functions
$viz.$, survival function, hazard rate function, reversed hazard
rate function etc. corresponding to an order statistic coming from
different heterogeneous populations are not very explicit in nature
and hence similar treatment as above cannot be used to handle these
problems. We are planning to study different ordering results for
the $k$-out-of-$n$ system of heterogeneous populations under
multiple-outlier models, and then extend these results to the
general model.

    %Further, there are different relative orderings available in the literature (cf. Hazra and Nanda (2016)). A similar study for all those orderings are under study and will be reported separately.}

\section*{Acknowledgements:} The authors are thankful to the Editor-in-Chief, the Associate Editor and the anonymous Reviewers for valuable suggestions which lead to an improved version of the manuscript. %The
%financial support from NBHM, Govt. of India (vide Ref. No.2/48(25)/2014/NBHM(R.P.)/R\&D II/1393 dt. Feb. 3, 2015) is dulyacknowledged by Asok K. Nanda.

\end{document}